\documentclass[11pt, a4paper]{article}

\usepackage[T1]{fontenc}
\usepackage[utf8]{inputenc} 
\usepackage[usenames]{color} 
\usepackage{graphicx} 
\usepackage{times} 
\usepackage{amsmath} 
\usepackage{amssymb}  
\usepackage{subcaption}

\usepackage{url}
\usepackage{multirow}
\usepackage{cite}

\usepackage{algorithm}
\usepackage{pgf}
\usepackage{pgfplots}
\usepgfplotslibrary{fillbetween}
\usepackage{tikz}
\usetikzlibrary{arrows,automata}
\usetikzlibrary{intersections}
\usetikzlibrary{calc}
\usetikzlibrary{shapes}
\usetikzlibrary{positioning}
\usetikzlibrary{decorations.text}
\pgfdeclarelayer{background}
\pgfdeclarelayer{foreground}
\pgfsetlayers{background,main,foreground}
\usetikzlibrary{calc,patterns,decorations.pathmorphing,decorations.markings}

\definecolor{morange}{rgb}{0.8500,0.3250,0.0980}
\definecolor{myellow}{rgb}{0.9290,0.6940,0.1250}
\definecolor{myblue}{rgb}{0,0.4392,0.7529}
\definecolor{mycol}{RGB}{19,48,128}
\usepackage{hyperref} 
\hypersetup{
	colorlinks = true,
	citecolor = mycol,
	linkcolor = mycol,,
	urlcolor = black
}

\usepackage{enumitem}

\newtheorem{theorem}{Theorem}

\newtheorem{cond}{Condition}

\newcommand{\rr}{\mathbb{R}}
\newcommand{\st}{\mathop{\rm s.t.}\nolimits}

\newcommand{\beqar}{\begin{eqnarray}}
\newcommand{\eeqar}{\end{eqnarray}}
\newcommand{\beqarno}{\begin{eqnarray*}}
\newcommand{\eeqarno}{\end{eqnarray*}}
\newcommand{\ba}[1]{\begin{array}{#1}}
\newcommand{\ea}{\end{array}}

\renewcommand{\vec}{\mathop{\rm vec}\nolimits}
\newcommand{\Var}{\mathop{\rm Var}\nolimits}

\newcommand{\smallmat}[1]{\left[ \begin{smallmatrix}#1 \end{smallmatrix} \right]}

\newcommand{\sps}{\hspace{0.3mm}}

\newcommand{\plantonly}{{{\tt plant only}}}

\newcommand{\plantonlyt}{{{\tt plant only}}}
\newcommand{\combined}{{\tt combined}} 
\newcommand{\combinedt}{{{\tt combined}}}

\makeatletter
\renewcommand{\@fnsymbol}[1]{}
\makeatother

\begin{document}

	\title{\LARGE Efficient identification of linear, parameter-varying, and nonlinear systems with noise models}
	\author{Alberto Bemporad$^*$\thanks{$^*$IMT School for Advanced Studies Lucca, Lucca, Italy. Email: \texttt{alberto.bemporad@imtlucca.it}} and Roland T\'oth$^\dag$\thanks{$^\dag$Control Systems Group, Eindhoven University of Technology, Eindhoven, The Netherlands, and Systems and Control Laboratory, HUN-REN Institute for Computer Science and Control, Budapest, Hungary. Email: \texttt{r.toth@tue.nl}}
\thanks{
This research was supported by the European Union within the framework of the National Laboratory for Autonomous Systems (RRF-2.3.1-21-2022-00002), by the Air Force Office of Scientific Research under award number FA8655-23-1-7061, and by the European Research Council (ERC), Advanced Research Grant COMPACT (Grant Agreement No. 101141351). Views and opinions expressed are however those of the authors only and do not necessarily reflect those of the European Union or the European Commission. Neither the European Union nor the granting authority can be held responsible for them. Corresponding author: A. Bemporad.
}
}


\maketitle

\thispagestyle{empty}    

\begin{abstract}
We present a general system identification procedure capable of estimating of a broad spectrum of state-space dynamical models, including \emph{linear time-invariant} (LTI), \emph{linear parameter-varying} (LPV), and nonlinear (NL) dynamics, along with rather general classes of noise models. Similar to the LTI case, we show that for this general class of model structures, including the NL case, the model dynamics can be separated into a \emph{deterministic process} and a \emph{stochastic noise} part, allowing to seamlessly tune the complexity of the combined model both in terms of nonlinearity and noise modeling. We parameterize the involved nonlinear functional relations by means of \emph{artificial neural-networks} (ANNs), although alternative parametric nonlinear mappings can also be used. To estimate the resulting model structures, we optimize a prediction-error-based criterion using an efficient combination of  
a constrained quasi-Newton approach and automatic differentiation, achieving training times in the order of seconds compared to existing state-of-the-art ANN methods which may require hours for models of similar complexity. We formally establish the consistency guarantees for the proposed approach and demonstrate its superior estimation accuracy and computational efficiency on several benchmark LTI, LPV, and NL system identification problems.
\end{abstract}


\section{Introduction}

\emph{System identification}, i.e., data-driven modelling, has been an important field of research in control engineering and in statistics, providing reliable approaches for obtaining models of poorly understood systems, filling the gaps where first-principle knowledge is lacking, and providing the machinery for adaptation and data-driven design of controllers \cite{8897147}. However, if one goes beyond the well-established \emph{linear time-invariant} (LTI) system
identification framework \cite{Ljung:1999}, considering, e.g., \emph{linear time-varying} (LTV) \cite{Niedz2000}, \emph{linear parameter-varying} (LPV) \cite{Toth2010SpringerBook}, or various \emph{nonlinear} (NL) models \cite{Giri2010}, a core question of the resulting identification problems is how to 
efficiently estimate from data the functional relations involved in these models.

While classical {\it beyond-LTI} methods used fixed basis-function parameterizations to express nonlinearities, with the recent progress of computational capabilities and machine learning/deep learning methods, a wide-range 
of learning approaches has been introduced, from \emph{Gaussian process} (GP) regression~\cite{Pillonetto2014} and \emph{support vector machines}~\cite{Suykens2002} to \emph{artificial neural networks} (ANN) based training (see ~\cite{PILLONETTO2025111907} for a recent overview), including \emph{state-space} (SS) recurrent networks such as LSTMs \cite{GONZALEZ2018485,LJUNG20201175}, SS-ANNs~\cite{Suykens96,RIBEIRO2020109158,MASTI2021109666,forgione2020dynonet,Toth22a,Bem25,Bem23}, etc. In particular, the latter class of methods successfully combines deep learning with state-space model identification concepts, achieving extreme accuracy levels on well established benchmarks, see \cite{Toth22a}. However, in the quest of solving the underlying function estimation problem with improved accuracy and reliability, many aspects of the classical identification theory have been sacrificed. While deep-learning methods can achieve remarkable results, the relating training time when no state-measurements are available is often too high compared to classical identification methods, requiring tens of hours, sometimes days, with the most sophisticated stochastic gradient methods, like Adam~\cite{KB14}, to train mid-sized models on training sequences of reasonable length. Using such methods in an identification cycle, which requires iterations on model structure selection, hyperparameter optimization on validation data, and experiment design, is simply impractical. 

Next to the computational challenges, many of the developed machine learning-based identification methods do not consider measurement or process noise affecting the estimation problem, or are based on simplistic noise settings / isolated scenarios, with a few exceptions, e.g.,~\cite{Toth22a,RIBEIRO2020109158}. The general understanding of process and noise models, and the related concept of one/$n$-step-ahead predictors and connected modelling choices, model structures, and the overall estimation concept,
which are well-developed for the LTI case, is largely missing in the general NL context. Besides theoretical
limitations, in practice it is not clear to the user how to co-estimate linear and nonlinear elements, select a combination of process and noise models, and scale up the identification process from linear to nonlinear modeling within the same framework.

\subsection{Contribution}
To overcome the aforementioned challenges, our first contribution is to provide an extension of the classical process and noise modelling framework of \cite{Ljung:1999} to the general NL case, showing that even for general nonlinear stochastic state-space models with innovation noise, the dynamics can be separated into a deterministic process and a stochastic noise model. This allows us to unify a spectrum of model-structure representations from LTI to NL classes, including externally-scheduled and self-scheduled LPV models, enabling to tune seamlessly the complexity of the model, with model structure selection both in terms of nonlinearity and noise modeling. We show that, in the posed nonlinear system identification setting, the inverse of the noise model can always be analytically computed, and that it is beneficial to learn this inverse noise model rather than its direct counterpart. 

Our second contribution is to extend the identification methodology recently introduced in~\cite{Bem25}, which is based
on L-BFGS-B optimization~\cite{BLNZ95} and the efficient auto-differentiation framework of JAX~\cite{JAX}, to identify the resulting combination of process and noise models under various LTI, LPV, and NL model classes and ANN parameterizations of the involved nonlinear functional relations.
As we demonstrate on various simulation and benchmark examples, the resulting identification approach 
is a highly efficient and unified approach to system identification that can achieve, and even surpass, the accuracy of the state-of-the-art deep-learning based identification methods, like SUBNET~\cite{Toth22a}, while bringing the total training time from the order of hours to the order of seconds. 

Finally, our last contribution is to show consistency of the resulting identification method. Together with formally proving that the dynamics of a nonlinear data-generating system with an innovation type of noise structure can be generally separated to a process model and nonlinear noise model that can be jointly estimated, this makes our approach a statistically and computationally reliable tool for solving identification problems in practice with the same level of simplicity of the existing tools available for the LTI case. 

The paper is structured as follows. In Section \ref{sec:2}, we introduce the considered identification setting and show our separation result. We also discuss the unified parameterization of model structures for capturing the dynamics 
of the data-generating system, and introduce the considered identification problem. This is followed by detailing the proposed identification method for the introduced model structures in Section \ref{sec:3}, together with a formal proof of its consistency in Section \ref{sec:4}. The estimation accuracy and computationally efficiency of the proposed identification approach is demonstrated in Section \ref{sec:examples} on several examples and benchmarks. Finally, conclusions about the achieved results are given in Section \ref{sec:6}.  

\subsection{Notation}
$\mathbb{R}$ and $\mathbb{Z}$ denote the set of real and integer numbers, respectively, while $\mathbb{Z}^+$ (or $\mathbb{Z}_0^+$) is the set of positive (non-negative) integers. Given a matrix $A\in \mathbb{R}^{n\times m}$, we denote by $\vec(A)$ the $nm$-dimensional column vector obtained by placing all the columns of $A$ below each other, and, more generally, given $p$ vectors $v_1,\ldots,v_p$ by $\vec(v_1,\ldots,v_p)=[v_1'\ \ldots\ v_m']'$. The set $\mathbb{I}_{\tau_1}^{\tau_2}$ denotes the set of integers $\tau$ such that $\tau_1\leq \tau \leq \tau_2$.

\section{Identification Setting} \label{sec:2}

\subsection{Data-generating system}
Consider the to-be-identified \emph{discrete-time} system defined in terms of a \emph{state-space} (SS) representation with an \emph{innovation}-type of noise process:
\begin{subequations} \label{model:innovation}
\begin{align}
w_{k+1} &= f_\mathrm{w}(w_{k},u_k,e_k),  \label{eq:system:proc:inv}\\
y_k & = g_\mathrm{w}(w_{k},u_k)+e_k, \label{eq:system:noise:inv}
\end{align}
\end{subequations}
where $k\in\mathbb{Z}$ denotes the discrete-time step, $u_k\in\rr^{n_u}$ a quasi-stationary input signal, $w_k\in\rr^{n_\mathrm{w}}$ the state, $y_k\in\rr^{n_y}$ the output, and $e_k$ is an i.i.d.~(identically and independently distributed) white noise process with finite variance $\Sigma_\mathrm{e}\in \mathbb{R}^{n_\mathrm{y} \times n_\mathrm{y}}$, independent of $u$. The functions $f_\mathrm{w}:\rr^{n_\mathrm{w}}\times\rr^{n_\mathrm{u}}\to\rr^{n_\mathrm{w}}$ and
$g_\mathrm{w}:\rr^{n_\mathrm{w}}\times\rr^{n_\mathrm{u}}\to\rr^{n_\mathrm{y}}$, i.e., the state transition and output functions, are considered to be real-valued,
deterministic, and possibly nonlinear functions. By assuming various structures for $f_\mathrm{w}$ and $g_\mathrm{w}$, many well-known model structures can be obtained such as \emph{nonlinear output noise} (NOE), \emph{nonlinear auto-regressive with exogenous input} (NARX), \emph{nonlinear auto-regressive with moving average exogenous input} (NARMAX) and \emph{nonlinear Box-Jenkins} (NBJ)~\cite{Billings2013NARMAX, jansson2003ARXsubspace}. For instance, if  $f_\mathrm{w}$ does not depend on $e_k$, then a \emph{nonlinear} (NL)-SS model with an OE noise structure is obtained.

Under mild conditions, the dynamics of \eqref{model:innovation} can be separated into a \emph{deterministic process} and a \emph{stochastic noise} part. Let
\begin{subequations}
\begin{equation}
G_\mathrm{o} : \left\{
\begin{aligned}
x_{k+1} &= f_\mathrm{x}(x_k,u_k),\\
 y_{\mathrm{o},k} &= g_\mathrm{x}(x_k,u_k),
\end{aligned} \right.
\label{eq:system:proc}
\end{equation}
correspond to the deterministic process part of~\eqref{model:innovation}, with $x_k\in\rr^{n_\mathrm{x}}$  being the associated state and $f_\mathrm{x}:\rr^{n_\mathrm{x}}\times\rr^{n_\mathrm{u}}\to\rr^{n_\mathrm{x}}$,
$g_\mathrm{x}:\rr^{n_\mathrm{x}}\times\rr^{n_\mathrm{u}}\to\rr^{n_y}$ being deterministic maps. To ensure an equivalence with \eqref{eq:system:noise:inv}, the output  $y_{\mathrm{o},k}$ is observed under a stochastic disturbance $v_k$:
\begin{equation}
y_k=y_{\mathrm{o},k} + v_k, \label{eq:system:output}
\end{equation} 
where $v_k$ is defined by the stable hidden Markov model:
\begin{equation}
H_\mathrm{o}: \left\{
\begin{aligned}
z_{k+1} &= f_\mathrm{z}(z_k,{x_k,u_k},e_k),\\
v_k &= g_\mathrm{z}(z_k,{x_k,u_k})+e_k.
\end{aligned} \right.
\label{eq:system:noise}
\end{equation}
\end{subequations}
In~\eqref{eq:system:noise}, $z_k$ is a latent stochastic process taking values in $\rr^{n_\mathrm{z}}$  at each time instant $k$, while $f_\mathrm{z}:\rr^{n_\mathrm{z}}\times \rr^{n_\mathrm{x}} \times \rr^{n_\mathrm{u}}  \times\rr^{n_\mathrm{y}}\to\rr^{n_\mathrm{z}}$ and $g_\mathrm{z}:\rr^{n_\mathrm{z}}{\times \rr^{n_\mathrm{x}} \times \rr^{n_\mathrm{u}} }\to\rr^{n_\mathrm{y}}$ are deterministic and possibly nonlinear functions. In the sequel, we will also treat $G_\mathrm{o}$ and $H_\mathrm{o}$ as dynamic operators, expressing the signal relations of \eqref{eq:system:proc} and \eqref{eq:system:noise} under left-compact support, i.e., we can write $y_{\mathrm{o},k}=G_\mathrm{o} u_k$ and $v_{k}=H_\mathrm{o}(x_k,u_k) e_k$. For the sake of simplicity of notation, we will omit the dependency of $H_\mathrm{o}$ on $(x_k,u_k)$ since now on.

\begin{theorem}[Separation of process and noise models]\label{thm:separation} A system given by  \eqref{model:innovation} always has a representation in the form of \eqref{eq:system:proc}--\eqref{eq:system:noise} with $n_\mathrm{x},n_\mathrm{z}\leq n_\mathrm{w}$, such that for any sample path realization of $e\in\mathbb{R}^{\mathbb{Z}_0^+}$ and $u\in\mathbb{R}^{\mathbb{Z}_0^+}$ and any initial condition $w_0\in\mathbb{R}^{n_\mathrm{w}}$, there exist $x_0\in\mathbb{R}^{n_\mathrm{x}}$ and $z_0\in\mathbb{R}^{n_\mathrm{z}}$ such that the resulting output trajectories in terms of \eqref{eq:system:output} and \eqref{eq:system:noise:inv} are equivalent.
\end{theorem} 

{\it Proof.} See Appendix \ref{proof:separation}. 

If additionally, in \eqref{eq:system:noise}, $f_\mathrm{z}(0,\centerdot,\centerdot,0)\equiv 0$ and $g_\mathrm{z}(0,\centerdot,\centerdot)\equiv 0$ are satisfied, then we call $G_\mathrm{o}$ and $H_\mathrm{o}$ a well-posed separation of the deterministic and stochastic parts of the system. See the proof of Theorem \ref{thm:separation} for such a construction. In the remaining part of the paper and in our examples, we will take this as an assumption without any loss of generality. We also highlight that, due to the nonlinearity of \eqref{eq:system:noise}, $\mathbb{E}_{e,z_0} \{ v_k\}=\bar{v}_k$ along a given $\{x_\tau\}_{\tau=0}^{k}$ and $\{u_\tau\}_{\tau=0}^{k}$ is often nonzero, where the expectation is taken  w.r.t. $e$ and an initial probability density function $p_\mathrm{z,0}$ of $z_0$. 

\begin{figure}
\centerline{ \includegraphics[scale=0.45]{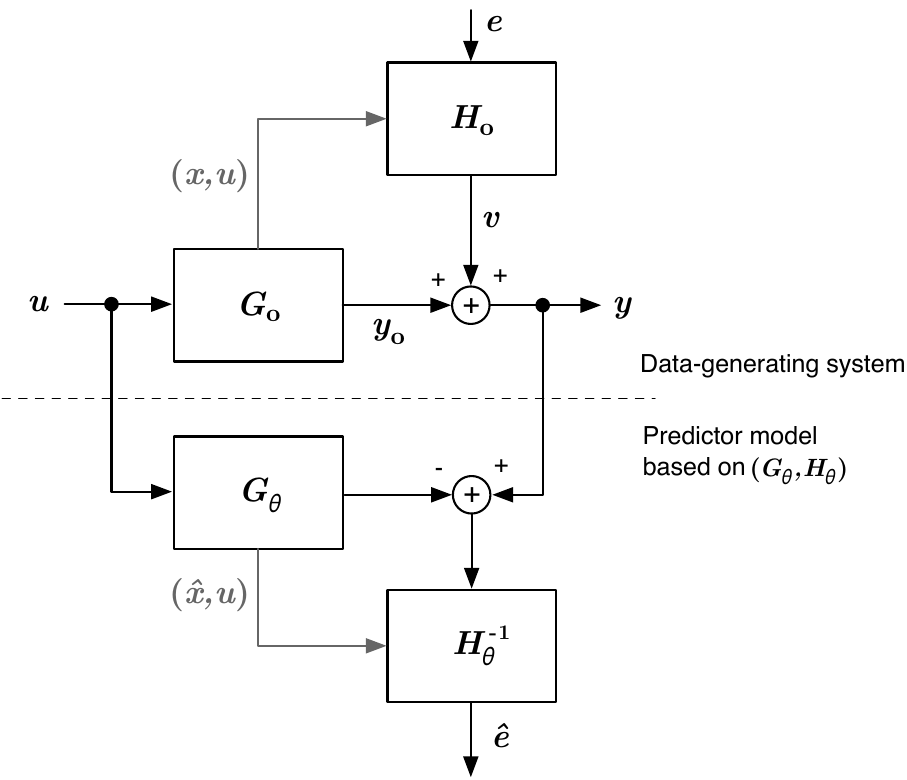}}
\caption{Separation of the data-generating system to a deterministic and stochastic part in the considered prediction-error-minimization setting.}\label{fig:pem}
\end{figure}

Equations \eqref{eq:system:proc}--\eqref{eq:system:noise} define the \emph{data-generating} system depicted in Figure \ref{fig:pem}, that we intend to identify from an \emph{input-output} (IO) dataset $\mathcal{D}_{N}=\{(u_k,y_k)\}_{k=0}^{N-1}$, where $\{ y_k \}_{k=0}^{N-1}$ is the response collected from the data-generating system~\eqref{eq:system:proc}--\eqref{eq:system:noise} in terms of a sample path realization for a given excitation sequence $\{ u_k \}_{k=0}^{N-1}$ and potentially unknown initial states $x_0\in\mathbb{R}^{n_\mathrm{x}}$ and $z_0\in\mathbb{R}^{n_\mathrm{x}}$. To avoid unnecessary clutter, we will not use different notation for random variables such as $e_k$ and their sampled values, but at places where confusion might arise, we will specify which notion is used.

As $e_k =v_k-g_\mathrm{z}(z_k, {x_k,u_k}),$ the inverse of the noise process \eqref{eq:system:noise}, can be expressed analytically as
\begin{equation}
H_\mathrm{o}^{-1}: \left\{ \begin{aligned}
z_{k+1} &= f_\mathrm{z}(z_k,{x_k,u_k},v_k-g_\mathrm{z}(z_k,{x_k,u_k})),\\
e_k &= v_k-g_\mathrm{z}(z_k,{x_k,u_k}).
\end{aligned} \right.
\label{eq:system:invnoise}
\end{equation}
It is easy to show that, for any $z_0\in\rr^{n_\mathrm{z}}$
and sequence $\{x_k,u_k\}_{k=0}^\infty$, we have
$e_k=H_\mathrm{o}^{-1} H_\mathrm{o} e_k$ and $v_k= H_\mathrm{o} H_\mathrm{o}^{-1} v_k$ 
when~\eqref{eq:system:noise} and~\eqref{eq:system:invnoise} are both initialized at $z_0$.
Moreover, by defining two functions $\tilde f_\mathrm{z}$ and $\tilde{g}_\mathrm{z}$ such that
$\tilde{f}_\mathrm{z}(z_k,{x_k,u_k},v_k) = f_\mathrm{z}(z_k,{x_k,u_k}, v_k-g_\mathrm{z}(z_k,{x_k,u_k}))$
and $g_\mathrm{z}=-\tilde{g}_\mathrm{z}$, we obtain
\[ 
    f_\mathrm{z}(z_k,{x_k,u_k},e_k)= \tilde{f}_\mathrm{z}(z_k,{x_k,u_k},e_k-\tilde g_\mathrm{z}(z_k,{x_k,u_k})).
\] By exploiting that $v_k=y_k- y_{\mathrm{o},k}= y_k -G_\mathrm{o}u_k $, we get
\begin{equation} \label{eq:noise:inv}
e_k=H_\mathrm{o}^{-1} (y_k-G_\mathrm{o}u_k).
\end{equation}
Next, using Equations \eqref{eq:system:proc}--\eqref{eq:system:noise}, we can write 
\begin{subequations}
\begin{align}
y_k & =G_\mathrm{o} u_k + H_\mathrm{o} e_k,\\
& = G_\mathrm{o} u_k + (H_\mathrm{o}-I) e_k + e_k,
\end{align}
and, by substituting \eqref{eq:noise:inv}, we get
\begin{align}
y_k 
& = G_\mathrm{o} u_k + (H_\mathrm{o}-I) H_\mathrm{o}^{-1} (y_k-G_\mathrm{o}u_k) + e_k,\\
& = G_\mathrm{o} u_k + (I - H_\mathrm{o}^{-1}) (y_k-G_\mathrm{o}u_k) + e_k, \label{eq:onestep:preform}
\end{align}
\end{subequations}
which remarkably holds true, as in the LTI case, despite of the possible nonlinearity of $G_\mathrm{o}$ and $H_\mathrm{o}$. Then, following the classical reasoning, under the information set of observed $y^{(k-1)}=\{y_\tau\}_{\tau\leq k-1}$ and $u^{(k)}= \{u_\tau\}_{\tau\leq k}$, the one-step-ahead prediction of $y(k)$ w.r.t.\ the $\ell_2$ loss is $\hat{y}(k\!\mid\! k-1)= \arg\ \min_{\delta\in\mathbb{R}^{n_\mathrm{y}}}\ \mathbb{E}\{ \| y(k)-\delta \|_2^2\mid u^{(k)},\ y^{(k-1)}\}$. Based on \eqref{eq:onestep:preform}, the straightforward application of the expectation operator yields
\begin{equation} \label{eq:1-step-ahead}
\hat{y}(k\!\mid\! k-1)= G_\mathrm{o} u_k + (I - H_\mathrm{o}^{-1}) (y_k-G_\mathrm{o}u_k) ,
\end{equation}  
which can be given in terms of an SS representation as
\begin{subequations}
\begin{align}
x_{k+1} &= f_\mathrm{x}(x_k,u_k), \\
z_{k+1} &= \tilde{f}_\mathrm{z}(z_k,{x_k,u_k},y_k-g_\mathrm{x}(x_k,u_k)), \\
\hat{y}(k\!\mid\! k-1) & = g_\mathrm{x}(x_k,u_k) - \tilde{g}_\mathrm{z}(z_k,{x_k,u_k}).
\end{align}
\end{subequations}
To the best of the authors' knowledge, such a general and computable form of the one-step-ahed predictor in the nonlinear case has not been shown in the literature before. 

\subsection{Model parameterization} \label{sec:statement}
Given the dataset $\mathcal{D}_N$, we aim to estimate $G_\mathrm{o}$ and $H_\mathrm{o}$ in terms of parametrized forms of \eqref{eq:system:proc} and \eqref{eq:system:noise}.
Specifically, we consider the model $G_\theta$ of the process in the following parametric SS form
\begin{subequations}
\begin{equation}
G_\theta : \left\{
\begin{aligned}
\hat{x}_{k+1} &= f_\mathrm{x}(\hat{x}_k,u_k,\theta_\mathrm{x}),\\
\hat y_k &= g_\mathrm{x}(\hat{x}_k,u_k,\theta_\mathrm{y}),
\end{aligned} \right.
\label{eq:model-x}
\end{equation}
where $\hat{x}_k\in\rr^{\hat{n}_\mathrm{x}}$ is the process model state, with $\hat{n}_\mathrm{x}\geq 1$ not necessarily being equal to (the usually unknown) number $n_\mathrm{x}$
of states of the data-generating system, and $f_\mathrm{x}:\rr^{\hat{n}_\mathrm{x}}\times\rr^{n_\mathrm{u}}\times\rr^{n_{\theta_\mathrm{x}}}\to\rr^{\hat{n}_\mathrm{x}}$,
$g_\mathrm{x}:\rr^{\hat{n}_\mathrm{x}}\times\rr^{n_\mathrm{u}}\times\rr^{n_{\theta_\mathrm{y}}}\to\rr^{n_\mathrm{y}}$ are functions parametrized by vectors $\theta_\mathrm{x}$ and $\theta_\mathrm{y}$, respectively (with a slight abuse of notation, for more clarity we use the same function names used to characterize the data-generating system).
We will specify later different parameterizations and structures of these functions. 

To deal with the presence of noise, i.e., to model the discrepancy 
\[
    \hat v_k=y_k-\hat y_k,
\]
with $\hat v_k\in\rr^{n_\mathrm{y}}$, between the measured output $y_k$ and the output $\hat y_k$ predicted by the process model, we aim to estimate the inverse of $H_\mathrm{o}$ in terms of the parameterized noise model  $H_\theta^{-1}$, described by
\begin{equation}
H_\theta^{-1} : \left\{
\begin{aligned}
\hat z_{k+1} &= \tilde f_\mathrm{z}(\hat z_k,{\hat{x}_k,u_k},\hat v_k,\theta_\mathrm{z}),\\
\hat e_k &=\tilde g_\mathrm{z}(\hat z_k,{\hat{x}_k,u_k},\theta_\mathrm{e})+\hat v_k,
\end{aligned} \right.
\label{eq:model-z-inv}
\end{equation}
where $\hat z_k\in\rr^{\hat{n}_\mathrm{z}}$ is the state of the inverse noise model, $\tilde f_\mathrm{z}:\rr^{\hat{n}_\mathrm{z}}\times \rr^{\hat{n}_\mathrm{x}} \times \rr^{n_\mathrm{u}} \times\rr^{n_\mathrm{y}}\times\rr^{n_{\theta_\mathrm{z}}}\to\rr^{n_\mathrm{z}}$ and $\tilde g_\mathrm{z}:\rr^{\hat{n}_\mathrm{z}}\times \rr^{\hat{n}_\mathrm{x}} \times \rr^{n_\mathrm{u}} \times\rr^{n_{\theta_\mathrm{e}}}\to\rr^{n_\mathrm{y}}$
are parametrized by two more vectors $\theta_\mathrm{z}$ and $\theta_\mathrm{e}$, respectively,
to be determined together with $\theta_\mathrm{x}$ and $\theta_\mathrm{y}$. 
Note that, based on \eqref{eq:system:invnoise}, the chosen model structure~\eqref{eq:model-z-inv}
enables us to easily derive the forward noise model $H_\theta$ as
\begin{equation}\hspace{-1mm}H_\theta : \left\{
\begin{aligned}
\hat z_{k+1} &\!=\! \tilde f_\mathrm{z}(\hat z_k,{\hat{x}_k,u_k}, \hat e_k\!-\!\tilde g_\mathrm{z}(z_k,{\hat{x}_k,u_k},\theta_\mathrm{e}),\theta_\mathrm{z}),\\
\hat v_k &\!=\! \hat e_k\!-\!\tilde g_\mathrm{z}(z_k,{\hat{x}_k,u_k}, \theta_\mathrm{e}),
\end{aligned} \right.
\label{eq:model-z}
\end{equation}
\label{eq:general-parametric-model}%
\end{subequations}
which gives that  $\tilde f_\mathrm{z}$ and $\tilde g_\mathrm{z}$ must satisfy the condition $\tilde g_\mathrm{z}(0,\centerdot,\centerdot,\theta_\mathrm{e})\equiv 0$ and $\tilde f_\mathrm{z}(0,\centerdot,\centerdot,0,\theta_\mathrm{z})\equiv 0$ for all $\theta_\mathrm{z}\in \Theta_\mathrm{z} \subseteq \rr^{n_{\theta_\mathrm{z}}} $ to satisfy the separation constraints. 
Different ways of parameterizing $G_\theta$ and $H_\theta^{-1}$ will be detailed in Section~\ref{sec:models}.

\subsection{Identification problem}
Based on the given dataset $\mathcal{D}_N$, the objective is to estimate the model parameters $\theta=\mathrm{vec}(\theta_\mathrm{x},\theta_\mathrm{y},\theta_\mathrm{z},\theta_\mathrm{e})$ by minimizing the criterion
\begin{equation} \label{eq:cost}
V_{\mathcal{D}_N}(\theta,\hat{w}_0)=  \frac{1}{N} \sum_{k=0}^{N-1} \|\hat{e}^\mathrm{pred}_k\|_2^2
\end{equation}
where $\hat{e}^\mathrm{pred}_k=y_k-\hat{y}^\mathrm{pred}_k$ and $\hat{y}^\mathrm{pred}_k$ is computed via the prediction model
\begin{subequations}
 \begin{align}
\hat{x}_{k+1} &= f_\mathrm{x}(\hat{x}_k,u_k,\theta_\mathrm{x}),\\
\hat{z}_{k+1} &= \tilde f_\mathrm{z}(\hat{z}_k,\hat{x}_k,u_k,y_k- g_\mathrm{x}(\hat{x}_k,u_k,\theta_\mathrm{y}),\theta_\mathrm{z}),\\
\hat{y}^\mathrm{pred}_k &= g_\mathrm{x}(\hat{x}_k,u_k,\theta_\mathrm{y}) -  \tilde g_\mathrm{z}(\hat z_k, \hat{x}_k,u_k,\theta_\mathrm{e}),
\end{align} 
\label{eq:1-step-prediction}%
\end{subequations}
for $k\in \mathbb{I}_0^{N-1}$, with $\hat{w}_0=\mathrm{vec}(\hat{x}_0,\hat{z}_0)$ as the initial condition, $\hat{w}_0\in\mathbb{R}^{\hat{n}_\mathrm{w}}$, $\hat{n}_\mathrm{w}=\hat{n}_\mathrm{x}+\hat{n}_\mathrm{z}$. 

The objective~\eqref{eq:1-step-prediction} corresponds to the traditional \emph{prediction-error-minimization} (PEM) setting and aims to minimize the error of the one-step-ahead predictor of the combined model structure in terms of \eqref{eq:1-step-ahead} w.r.t. the observed output samples $\{y_k\}_{k=0}^{N-1}$.

\section{Model structures}
\label{sec:models}
Next, we discuss various parameterizations of $G_\theta$ and $H_\theta^{-1}$ that can be used to estimate linear, parameter-varying, and nonlinear models.

\subsection{LTI models}
To restrict the model structure to only express LTI dynamics, we take $G_\theta$ in the following form \vspace{-2mm}
\begin{subequations} \begin{equation}
G_\theta : \Biggl\{
\begin{aligned}
\hat{x}_{k+1}& =\overbrace{A_\mathrm{x} \hat{x}_k+B_\mathrm{x} u_k}^{ f_\mathrm{x}(\hat{x}_k,u_k,\theta_\mathrm{x}) }, \\
\hat{y}_{k} & = \underbrace{C_\mathrm{x} \hat  x_k+D_\mathrm{x} u_k}_{g_\mathrm{x}(\hat{x}_k,u_k,\theta_\mathrm{y})},
\end{aligned}
\label{eq:LTI-x} \vspace{-3mm}
\end{equation}
and parametrize $H_\theta^{-1}$ as \vspace{-2mm}
\begin{equation}
H_\theta^{-1} : \Biggl\{
\begin{aligned}
\hat z_{k+1} &= \overbrace{A_\mathrm{z} \hat z_k+B_\mathrm{z} \hat v_k,}^{\tilde f_\mathrm{z}(\hat z_k,{\hat{x}_k,u_k},\hat v_k,\theta_\mathrm{z})} \\
\hat e_k &= \hspace{-4mm} \underbrace{C_\mathrm{z}\hat z_k}_{\tilde g_\mathrm{z}(\hat z_k,{\hat{x}_k,u_k},\theta_\mathrm{e})}\hspace{-5mm} +\ \hat v_k,
\end{aligned} 
\label{eq:LTI-z-inv}
\end{equation} 
\label{eq:model-LTI}%
\end{subequations}
for which $\theta_\mathrm{x}=\vec([A_\mathrm{x}\ B_\mathrm{x}])$,  $\theta_y=\vec([C_\mathrm{x} \ D_\mathrm{x}])$, $\theta_\mathrm{z}=\vec([A_\mathrm{z}\ B_\mathrm{z}])$,
$\theta_\mathrm{e}=\vec(C_\mathrm{z})$. Note that, by  appropriately restricting $(A_\mathrm{z},B_\mathrm{z},C_\mathrm{z})$ further, well known LTI model structures 
can be realized, from ARX to BJ\footnote{For ARX and ARMAX models, an observer canonical realization form can be used to separate the states associated with the autoregressive and moving average dynamics, see \cite[Ch. 4]{Ljung:1999}; OE models require setting $(A_\mathrm{z},B_\mathrm{z},C_\mathrm{z}) = 0$; BJ corresponds to an unrestricted parametrization of $(A_\mathrm{z},B_\mathrm{z},C_\mathrm{z})$.}. 

\subsection{LPV models}
To estimate the system dynamics in LPV form, $G_\theta$ and $H_\theta^{-1}$ are formulated as
\begin{subequations} \begin{equation}
G_\theta : \Biggl\{
\begin{aligned}
\hat{x}_{k+1}& = A_\mathrm{x}(p_k) \hat{x}_k+B_\mathrm{x}(p_k) u_k, 
\\
\hat{y}_{k} & = C_\mathrm{x}(p_k) \hat  x_k+D_\mathrm{x}(p_k)u_k,  
\end{aligned}
\label{eq:LPV-x} \vspace{-1mm}
\end{equation}
\begin{equation}
\hspace{-4.5mm} H_\theta^{-1} : \Biggl\{
\begin{aligned}
\hat z_{k+1} &=A_\mathrm{z}(p_k) \hat z_k+B_\mathrm{z}(p_k) \hat v_k, 
\\
\hat e_k &= C_\mathrm{z}(p_k)\hat z_k + \hat v_k,
\end{aligned} 
\label{eq:LPV-z-inv}
\end{equation} \end{subequations}
where $p_k\in\rr^{n_\mathrm{p}}$ is the scheduling variable. In the so called \emph{self-scheduled} case (see \cite{Toth22CDCb} for a detailed discussion on LPV modelling and identification), $p_k=\psi(x_k,u_k,\theta_\psi)$, with $\psi:\rr^{n_\mathrm{x}}\times \rr^{n_\mathrm{u}} \times \rr^{n_{\theta_\psi}}\rightarrow \rr^{n_\mathrm{p}}$ being a parametrized, generally nonlinear function, called the \emph{scheduling map}. In the alternative \emph{externally-scheduled} case,
$p_k$ is considered as a given exogenous signal, which is part of the input $\tilde{u}_k=\vec(u_k,p_k)$ received by the system, and is included in the dataset $\mathcal{D}_N$. 
By introducing,
\[
    M_\mathrm{x}(p_k)=\smallmat{A_\mathrm{x}(p_k) & B_\mathrm{x}(p_k)\\C_\mathrm{x}(p_k) & D_\mathrm{x}(p_k)}\!,\ \ M_\mathrm{z}(p_k)=\smallmat{A_\mathrm{z}(p_k) & B_\mathrm{z}(p_k)\\C_\mathrm{z}(p_k) & 0}\!,
\]
we parameterize the LPV model as 
\begin{equation}
M_\mathrm{x}(p_k)\!=\!f_{M_\mathrm{x}}(p_k,\theta_{M_\mathrm{x}}),\  M_\mathrm{z}(p_k)\!=\!f_{M_\mathrm{z}}(p_k,\theta_{M_\mathrm{z}}), 
\label{eq:model-lpv-fnn}
\end{equation}%
where $f_{M_\mathrm{x}}$ and $f_{M_\mathrm{z}}$ are multi-layer \emph{feedforward neural networks} (FNNs) with $L-1\geq 0$ hidden layers and a linear output layer having  
$n_p$ inputs and $(n_\mathrm{x}+n_\mathrm{y})(n_\mathrm{x}+n_\mathrm{u})$ and $n_\mathrm{z}^2+2n_\mathrm{y}n_\mathrm{z}$ outputs, respectively, and
parameterized by the vectors $\theta_\mathrm{M_x}$ and $\theta_\mathrm{M_z}$, containing the weight and bias terms of the layers. Also, these FNNs can also contain a linear bypass for more efficient modeling of linearly dominant functions. Note that different nonlinear parametric functions can be used here in alternative to FNNs, 
including parametric physics-based models.

A special case of~\eqref{eq:model-lpv-fnn} when $f_{M_\mathrm{x}}$ and $f_{M_\mathrm{z}}$ only consist of a linear layer, i.e., $L=1$: 
\begin{equation}
M_\mathrm{x}(p_k)=M_{\mathrm{x},0}+\sum_{i=1}^{n_\mathrm{p}} p_{k,i}M_{\mathrm{x},i}
\label{eq:model-lpv}
\end{equation}%
and similarly for $M_\mathrm{z}(p_k)$, corresponding to LPV models that depend \emph{affinely} on $p_k$. The class of LPV state-space models with affine dependence has paramount importance in control engineering, due the wide range of convex analysis and controller synthesis methods that assume the representation of the system in such a form.

In the self-scheduled setting, the function $\psi$ defining $p_k$ is taken as a multi-layer FNN with weights and biases collected in $\theta_\psi$, while restricting $f_{M_\mathrm{x}}$ and $f_{M_\mathrm{z}}$
to the linear expression~\eqref{eq:model-lpv}. This is without loss of generality:
due to the flexibility we have in choosing function $\psi$, given an arbitrary FNN representation of
$f_{M_\mathrm{x}}$ and $f_{M_\mathrm{z}}$, we can interpret its last linear layer as in~\eqref{eq:model-lpv},
set $p_k$ as the value of the neurons at the second last nonlinear layer, and call $\psi$ the remaining layers. The parameter vector to learn is therefore $\theta= \vec([M_{\mathrm{x},0},\ldots,M_{\mathrm{x},n_\mathrm{p}}], [M_{\mathrm{z},0},\ldots,M_{\mathrm{z},n_\mathrm{p}}], \theta_\psi )$. Note that, with reference to the general
model structure~\eqref{eq:general-parametric-model}, in this case,
$\theta_\mathrm{x}$, $\theta_\mathrm{y}$, $\theta_\mathrm{v}$, and $\theta_\mathrm{e}$ share the
components of $\theta_\psi$.

In the externally scheduled case, $p_k$ is given, therefore there is no need to estimate $\psi$. Instead,  
we take $f_{M_\mathrm{x}}$ and $f_{M_\mathrm{z}}$ as general FNNs, with $\theta = \vec(\theta_\mathrm{M_x},\theta_\mathrm{M_z})$ as the overall parameter vector, and $\theta_\mathrm{x}$, $\theta_\mathrm{y}$, $\theta_\mathrm{v}$, $\theta_\mathrm{e}$ defined accordingly with possibly shared components. 

\subsection{General nonlinear models} 
In the general case of black-box nonlinear models, $G_\theta$ and $H_\theta^{-1}$ are considered in the form of  
\eqref{eq:model-x} and \eqref{eq:model-z-inv}, corresponding to  \emph{recurrent neural network} (RNN) models  when $f_\mathrm{x}$, $g_\mathrm{x}$, $\tilde f_\mathrm{z}$, $\tilde g_\mathrm{z}$ are taken as multi-layered FNNs with linear output layer and a possible linear bypass.
The weights and bias terms of these networks are collected in $\theta_\mathrm{x}$, $\theta_\mathrm{y}$, $\theta_\mathrm{z}$, and $\theta_\mathrm{e}$, respectively. 

Note that, while in the LTI and LPV cases, the considered model structures naturally satisfy the separation constraint, in the general NL case we need to enforce $\tilde g_\mathrm{z}(0,\centerdot,\centerdot,\theta_\mathrm{e})\equiv 0$ and $\tilde f_\mathrm{z}(0,\centerdot,\centerdot,0,\theta_\mathrm{z})\equiv 0$ for all $\theta_\mathrm{z}\in \Theta_\mathrm{z} \subseteq \rr^{n_{\theta_\mathrm{z}}}$. 

\section{Identification method} \label{sec:3}

Next, we discuss the proposed identification method to estimate a model according to the parametrized model structure \eqref{eq:model-x} and \eqref{eq:model-z-inv} based on a given data set $\mathcal{D}_N$ and under the considered criterion~\eqref{eq:cost}. We provide a general method that can be applied to any of the model classes described in Section~\ref{sec:models}.

\subsection{Combined plant and noise model learning} 
\label{sec:opt}
To address the posed identification problem in terms of the minimization of \eqref{eq:cost} with respect to $\theta\in\rr^{n_\theta}$, which collects all the free model parameters, 
and $\hat{w}_0$, which collects the initial states of the plant and noise models,
we solve the following optimization problem
\begin{align}
\min_{\theta,\hat{w}_0} &\ \ V_{\mathcal{D}_N}(\theta,\hat{w}_0) + R(\theta,\hat{w}_0), 
\label{eq:optim-prob}  \\
\st&\left\{\ba{rcl}
\hat{x}_{k+1} &=& f_\mathrm{x}(\hat{x}_k,u_k,\theta_\mathrm{x}), \\
\hat{z}_{k+1} &=& \tilde f_\mathrm{z}(\hat{z}_k,\hat{x}_k,u_k,y_k- g_\mathrm{x}(\hat{x}_k,u_k,\theta_\mathrm{y}),\theta_\mathrm{z}),\\
\hat{e}^\mathrm{pred}_k &=& y_k- g_\mathrm{x}(\hat{x}_k,u_k,\theta_\mathrm{y}) +  \tilde g_\mathrm{z}(\hat z_k, \hat{x}_k,u_k,\theta_\mathrm{e}), 
    \ea\right. \notag
\end{align} 
where $R(\theta,\hat{w}_0)$ is a regularization function.  The term
\begin{equation}
R(\theta,\hat{w}_0) = \frac{\rho_\theta}{2}\|\theta \|_2^2 + \tau \|\theta \|_1 
+ \frac{\rho_\mathrm{w}}{2}\|\hat{w}_0 \|_2^2
\label{eq:regularization}
\end{equation}
includes an \emph{elastic-net} type of regularization on the model parameters $\theta$, where $\rho_\theta>0$ and $\tau\geq 0$ are the weights defining the $\ell_2$ and $\ell_1$ penalties, respectively, and $\rho_\mathrm{w}>0$ relates to the $\ell_2$-regularization of the initial state $\hat{w}_0$.

In this paper, we adopt the optimization approach described in~\cite{Bem25} based on handling
$\ell_1$ penalties by splitting $\theta=\theta_+-\theta_-$, $\theta_+,\theta_-\geq 0$, and 
bound-constrained limited-memory Broyden-Fletcher-Goldfarb-Shanno (L-BFGS-B) optimization~\cite{BLNZ95} to solve~\eqref{eq:optim-prob}. To better escape local minima, a fixed number of Adam~\cite{KB14} 
gradient-descent steps is used to warm-start the L-BFGS-B optimization (see~\cite[Table 3]{Bem25} for the benefit of using Adam as a warm start). We use automatic differentiation to compute the gradients of the objective function with respect to $(\theta,\hat w_0)$. In Section~\ref{sec:examples}, we will demonstrate that this optimization toolchain, extended from~\cite{Bem25} for the proposed model classes, gives a computationally far more efficient training of the considered ANN-based state-space models than existing nonlinear systems identification methods using stochastic gradient methods, e.g., Adam, directly applied on various forms of ANN models.  

\subsection{Model structure selection} \label{sec:mod:select}
In~\eqref{eq:regularization}, we  can also include a \emph{group-lasso} penalty 
\begin{equation}
    \tau_\mathrm{g} \sum_{i=1}^{n_g}\|\theta_{\mathrm{g},i}\|_2 
\label{eq:group-lasso}
\end{equation}
where  $\tau_\mathrm{g}>0$ is the regularization weight and $\theta_{\mathrm{g},i}$ represents groups of parameters selected from $\theta$ and $\hat{w}_0$ that are related to a specific state dimension of $x$ or $z$, or elements of the scheduling variable $p$. For example, in case of a self-scheduled LPV model, $ \theta_{\mathrm{g},i}$ can collect the parameters associated with the $i^\mathrm{th}$  dimension of $x$, containing the $i^\mathrm{th}$ rows and columns of the matrices  $\{A_{\mathrm{x},j}\}_{j=0}^{n_\mathrm{p}}$, the $i^\mathrm{th}$ rows of $\{B_{\mathrm{x},j}\}_{j=0}^{n_\mathrm{p}}$, the $i^\mathrm{th}$ columns of $\{C_{\mathrm{x},j}\}_{j=0}^{n_\mathrm{p}}$ and the $i^\mathrm{th}$ column of the input layer weights of $\psi$ and the $i^\mathrm{th}$ element of $\hat{w}_0$.
The additional penalty~\eqref{eq:group-lasso} can be used, in principle, to provide automated model structure selection (appropriate choice of $\hat{n}_\mathrm{x}$, $\hat{n}_\mathrm{z}$, arguments of the scheduling map,
groups of neurons in the NN layers, etc.) if the corresponding parameters can be zeroed without increasing the prediction cost.
See also~\cite{Bem25,MB24} for further examples on group-lasso regularization applied to system identification problems of nonlinear state-space models.
 
A simple implementation of model-structure selection based on group-lasso regularization involves: $i$) train a model-parameter vector  $\hat{\theta}_\mathrm{g}$ with $\tau_\mathrm{g}>0$; $ii$) remove superfluous optimization variables by zeroing the entries in $\hat{\theta}_\mathrm{g}$ whose absolute value is below a certain threshold $\epsilon_\mathrm{g}>0$; $iii$) re-estimate the remaining model 
parameters with $\tau_\mathrm{g}=0$. The process is repeated with different choices of the hyperparameters $\tau_\mathrm{g}$ and $\epsilon_\mathrm{g}$ until an accurate enough, but low-complexity, model is found.
This process can be also automated by $n$-fold cross-validation. As an alternative, suggested in~\cite{Candes2008}, an efficient process to recover lower-dimensional models is to weight the parameters appearing in the group-lasso term with the inverse of the corresponding values estimated at the previous estimation step and keep iterating until convergence is achieved.

\subsection{Initialization}
\label{sec:mod:init}
Since problem~\eqref{eq:optim-prob} is nonconvex,  it is often beneficial to apply bootstrapping, 
first training a deterministic model as in~\eqref{eq:model-x} by minimizing $\|y_k-\hat y_k\|_2^2=\|\hat v_k\|_2^2$
to get a warm start for $\theta_\mathrm{x}$, $\theta_\mathrm{y}$, and
then the combined plant and noise models~\eqref{eq:model-x} and~\eqref{eq:model-z-inv} by minimizing~\eqref{eq:cost}. 

To train NN functions, we initialize the bias terms to zero and the weights with values drawn from the normal distribution with zero mean. For the latter, a slightly more advanced initialization in the form of the Xavier method \cite{glorot2010Xavierinit} can be also applied. In case of LTI and LPV models, the elements of the additional matrices can be drawn from a normal distribution, although other known initialization schemes for state-space model identification can be followed, e.g., by pre-estimating models via subspace identification in the LTI case. 

\section{Consistency analysis} \label{sec:4}

In this section, we will prove the consistency properties of the proposed identification scheme which  corresponds to the minimization of \eqref{eq:optim-prob} under no regularization, i.e., $R(\theta,\hat w_0)=0$. In fact, we will show that the resulting model estimates will converge to an equivalent representation of the system in the form \eqref{eq:system:proc}-\eqref{eq:system:proc} if the number of data points in the available data set $\mathcal{D}_N$ tends to infinity, i.e., $N\rightarrow \infty$. 

\subsection{Convergence}
First we show convergence of the model estimates relying on the results in \cite{ljung1978convergence}. For this, we require the data generating system to satisfy  the following stability condition:

\begin{cond}[Stability of the data-generating system]
\label{assum:S1}
    The data-generating system \eqref{eq:system:proc}-\eqref{eq:system:noise} has the property that, for any $\delta> 0$, there exist constants $c\in [0, \infty)$ and $\lambda \in[0,1)$ such that
    \begin{align}
        \mathbb{E}_e \{ \|y_k - \tilde{y}_k\|_2^4 \} < c  \lambda^{k-k_\mathrm{o}}, \quad \forall k\geq k_\mathrm{o},
    \end{align}
    under any $k_\mathrm{o}\geq0$,  $ w_\mathrm{o},\tilde{w}_\mathrm{o}\in\mathbb{R}^{n_\mathrm{x}+n_\mathrm{z}}$ with $\|w_\mathrm{o}-\tilde{w}_\mathrm{o}\|_2<\delta$, and $\{(u_\ell,e_\ell)\}_{\ell={0}}^k\in \mathcal{S}_{[0,\infty]}$, where $\mathcal{S}_{[0,\infty]}$ denotes the $\sigma$-algebra generated by the random variables $\{(u_\ell,e_\ell)\}_{\ell={0}}^\infty$, and 
    $y_k$ and $\tilde{y}_k$ satisfy \eqref{eq:system:proc}-\eqref{eq:system:noise}  with the same $(u_k,e_k)$, but with initial conditions $\mathrm{vec}(x_{k_\mathrm{o}},z_{k_\mathrm{o}})=w_\mathrm{o}$ and $\mathrm{vec}(\tilde x_{k_\mathrm{o}},\tilde{z}_{k_\mathrm{o}})=\tilde w_\mathrm{o}$.
\end{cond}

To identify  \eqref{eq:system:proc}-\eqref{eq:system:noise}, let the \emph{model structure} $M_\theta = (G_\theta, H_\theta)$, defined by \eqref{eq:model-x}-\eqref{eq:model-z}, be chosen in terms of the discussed cases in Section \ref{sec:models}, corresponding to LTI, LPV with external or internal scheduling, or fully NL model structures with NN parametrization of the underlying nonlinear functions. Let $\theta$ be restricted to a compact set $\Theta\subset \mathbb{R}^{n_\theta}$, giving the \emph{model set} $\mathcal{M}=\{ M_\theta \mid \theta\in\Theta\} $. For a given initial condition $\hat{w}_0 \in \mathbb{R}^{\hat{n}_\mathrm{w}}$, Eq.~\eqref{eq:1-step-prediction} defines the following \emph{1-step-ahead predictor} of $M_\theta$
\begin{equation}
\hat{y}^\mathrm{pred}_k = m^\mathrm{pred}_k(\{y_\ell\}_{\ell=0}^{k-1},\{u_\ell\}_{\ell=0}^k,\theta, \hat{w}_0). \label{eq:predic:comp}
\end{equation}
Without loss of generality, we will take the assumption that $m^\mathrm{pred}_k$ is differentiable w.r.t. $\theta$ and $\hat{w}_0$ everywhere on an open neighborhood  $\breve{\Theta}$ of $\Theta$ and $\mathbb{R}^{\hat{n}_\mathrm{w}}$, respectively. Note that this is a technical condition, as only such parameterizations are considered for which automated differentiation is available as discussed in Section \ref{sec:opt}. Furthermore, to conduct our analysis, we require \eqref{eq:predic:comp} to be stable w.r.t. perturbations of the dataset, as formalized by the following condition:
\begin{cond}[Predictor stability]
\label{assum:M1} Given a $\delta\geq 0$, there exist constants $c\in [0, \infty)$ and $\lambda \in[0,1)$ such that,
for any $\theta\in\breve\Theta$, $\hat{w}_0\in\mathbb{R}^{\hat{n}_\mathrm{w}}$ with $\|\hat w_0\|_2 < \delta$, and for any $\{(y_\ell,u_\ell)\}_{\ell=0}^k,\{(\tilde{y}_\ell,\tilde{u}_\ell)\}_{\ell=0}^k \in\mathbb{R}^{(n_\mathrm{y}+n_\mathrm{u})\times (k+1)}$,
the predictors
  \begin{align*} \hat{y}^\mathrm{pred}_k &= m^\mathrm{pred}_k(\{y_\ell\}_{\ell=0}^{k-1},\{u_\ell\}_{\ell=0}^k,\theta, \hat{w}_0), \\ 
  \tilde{y}^\mathrm{pred}_k &= m^\mathrm{pred}_k(\{\tilde y_\ell\}_{\ell=0}^{k-1},\{\tilde u_\ell\}_{\ell=0}^k,\theta, \hat{w}_0),
  \end{align*} satisfy
    \begin{equation} \label{eq:M1}
        \|\hat{y}^\mathrm{pred}_k-\tilde{y}^\mathrm{pred}_k\|_2   \leq c \gamma(k) ,  \quad \forall k\geq 0,
\end{equation}
with $\gamma(k)=\sum_{\ell=0}^{k} \lambda^{k-\ell} \left(\|{u}_\ell-\tilde{u}_\ell\|_2 + \|{y}_\ell-\tilde{y}_\ell\|_2 \right)$
 and $\|m^\mathrm{pred}_k(\{0\}_{\ell=\hat{w}_0}^{k-1}$, $\{0\}_{\ell=0}^k$, $\theta$, $\hat{w}_0)\|_2 \leq c$.
Furthermore, \eqref{eq:M1} is also satisfied by $\frac{\partial}{\partial \theta}m^\mathrm{pred}_k$ and $\frac{\partial}{\partial \hat{w}_0}m^\mathrm{pred}_k$. 
\end{cond}

\begin{theorem}[Convergence] \label{lem:convergence}
Given the data-generating system \eqref{eq:system:proc}-\eqref{eq:system:noise}  satisfying Condition \ref{assum:S1} with a quasi-stationary $u$ independent of the white noise process ${e}$ and the model set $\mathcal{M}$ defined by  \eqref{eq:model-x}-\eqref{eq:model-z}, satisfying Condition \ref{assum:M1}. Then,
    \begin{align} \label{eq:conv}
        \underset{\theta\in \Theta, \hat{w}_0 \in \mathbb{R}^{\hat{n}_\mathrm{w}}}{\mathrm{sup}} \left \|V_{\mathcal{D}_N}(\theta,\hat{w}_0) - \mathbb{E}_e\{ V_{\mathcal{D}_N}(\theta,\hat{w}_0)\}\right\|_2 \rightarrow 0,
    \end{align}
    with probability 1 as $N \rightarrow \infty$. 
\end{theorem}

{\it Proof.} 
The identification criterion \eqref{eq:cost} satisfies  Condition C1 in \cite{ljung1978convergence}, hence the proof of \cite[Lemma 3.1]{ljung1978convergence} applies for the considered case.  \hfill $\blacksquare$

\subsection{Consistency}
In order to show consistency, we need to assume that the data-generating system $(G_\mathrm{o}, H_\mathrm{o})$ belongs to the chosen model set $(G_\theta, H_\theta)$. This means that for the considered model set $\mathcal{M}$ defined by  \eqref{eq:model-x}-\eqref{eq:model-z} under $\theta \in \Theta$, there exists a $\theta_\mathrm{o} \in \Theta$ such that the one-step-ahead predictor $m^\mathrm{pred}_k$ associated with $M_{\theta_\mathrm{o}}$ and the one-step-ahead predictor \eqref{eq:1-step-ahead} associated with the data-generating system $(G_\mathrm{o},H_\mathrm{o})$, defined by  \eqref{eq:system:proc}-\eqref{eq:system:noise}, correspond to equivalent state-space representations. However, as state-space representations are not unique, we can define $\Theta_\mathrm{o} \subset \Theta$ which contains all $\theta_\mathrm{o} \in \Theta$ that defines an equivalent model with the data-generating system $(G_\mathrm{o},H_\mathrm{o})$  in the one-step-ahead prediction sense. If $\Theta_\mathrm{o} = \varnothing$, then the model set is not rich enough to capture $(G_\mathrm{o},H_\mathrm{o})$ and consistency cannot be investigated.

We need the following condition to hold for the $\mathcal{D}_N$:
\begin{cond}[Persistence of excitation]
    \label{assum:PE} Given the model set $\mathcal{M}=\{ M_\theta \mid \theta\in\Theta\} $, we call the input sequence $\{u_\ell\}_{\ell=0}^{N-1}$ in $\mathcal{D}_{N}$ generated by $(G_\mathrm{o},H_\mathrm{o})$ {\it weakly persistently exciting}, if for all pairs of parameterizations given by $\left( \theta_1, \hat{w}_1 \right) \in \Theta \times \mathbb{R}^{\hat{n}_\mathrm{w}}$ and $\left( \theta_2, \hat{w}_2 \right) \in \Theta \times \mathbb{R}^{\hat{n}_\mathrm{w}}$ for which the function mapping is unequal, i.e., $V_{(\cdot)}(\theta_1,\hat{w}_1) \neq V_{(\cdot)}(\theta_2,\hat{w}_2)$, we have
    \begin{align}
        V_{\mathcal{D}_{N}}(\theta_1,\hat{w}_1) \neq V_{\mathcal{D}_{N}}(\theta_2,\hat{w}_2)
    \end{align}
    with probability 1.
\end{cond}

To show consistency, we also require that any element of $\Theta_\mathrm{o}$ and the corresponding initial condition $\hat w_\mathrm{o}$,  have minimal asymptotic cost in terms of $\lim_{N\rightarrow \infty}$ $V_{\mathcal{D}_N}(\theta,\hat{w}_0)$. Due to the mean squared loss function and since $\mathbb{E}_e \{ V_{\mathcal{D}_{N}}(\theta,\hat{w}_0) \}$ exists as shown in Theorem \ref{lem:convergence}, this trivially holds. For a detailed proof see \cite{ljung1978convergence}. 

\begin{theorem}[Consistency] \label{lem:consistency}
Under the conditions of Theorem \ref{lem:convergence} and Condition \ref{assum:PE}, 
    \begin{align}
        \underset{N \rightarrow \infty}{\lim} \hat{\theta}_N \in \Theta_\mathrm{o}
    \end{align}
    with probability 1, where
    \begin{align}
        (\hat{\theta}_N,\hat{w}_N) =\underset{\theta \in \Theta,{w}\in \mathbb{R}^{\hat{n}_\mathrm{w}}} { \arg \min} 
        V_{\mathcal{D}_{N}}(\theta, w).
    \end{align}
\end{theorem}
{\it Proof.}
The proof follows form Lemma 4.1 in \cite{ljung1978convergence} as the mean-squared loss function \eqref{eq:cost} fulfills Condition (4.4) in \cite{ljung1978convergence}. \hfill $\blacksquare$ 

\section{Examples} \label{sec:examples}
\begin{figure}[t]
\centering
\begin{subfigure}[t]{0.5\columnwidth} 
\centering
\includegraphics[height=2.5cm]{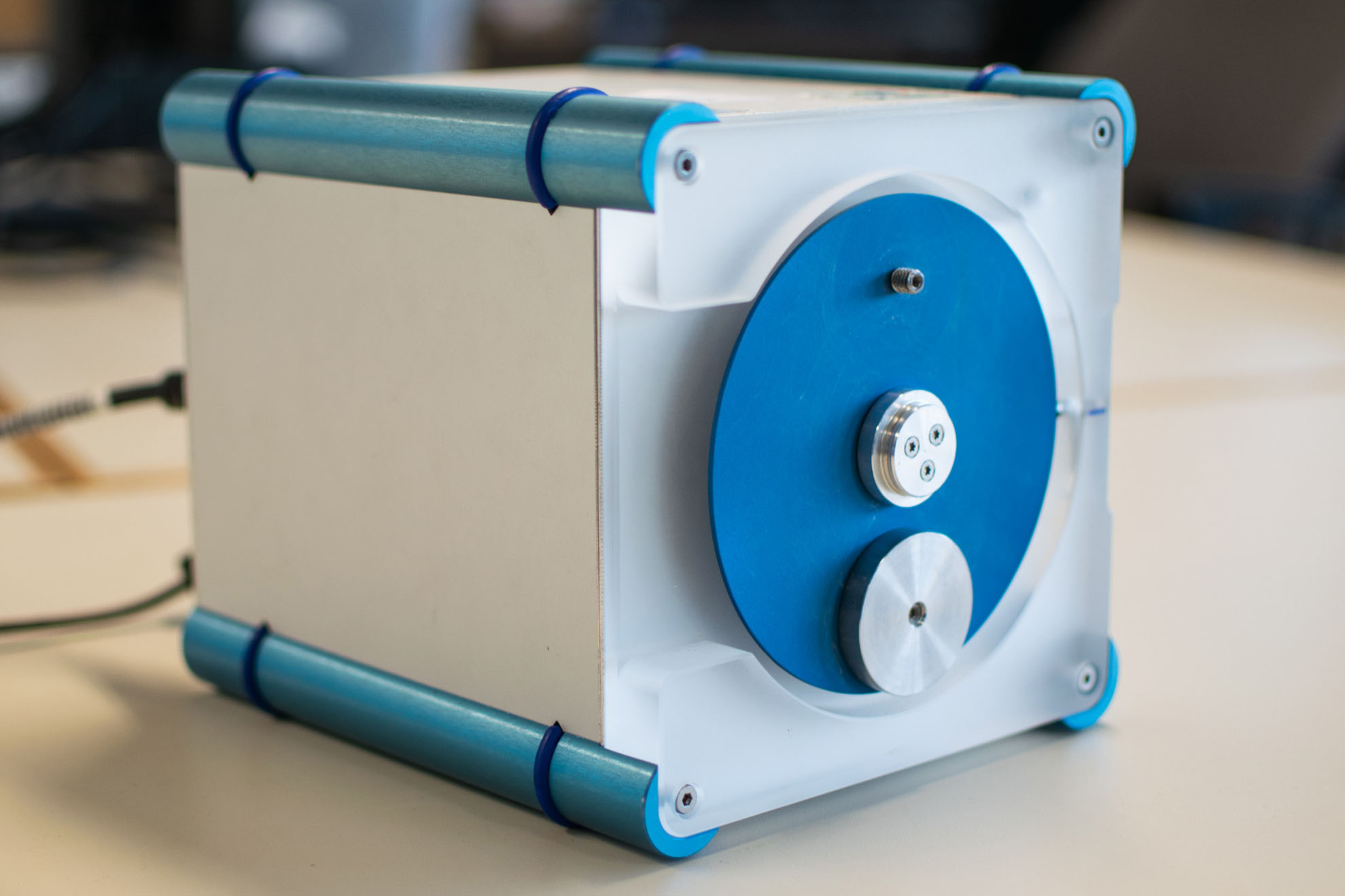}
\caption{Unbalanced disk} \label{fig:unbd}
\end{subfigure}\begin{subfigure}[t]{0.5\columnwidth}
\centering
\includegraphics[height=3.5cm]{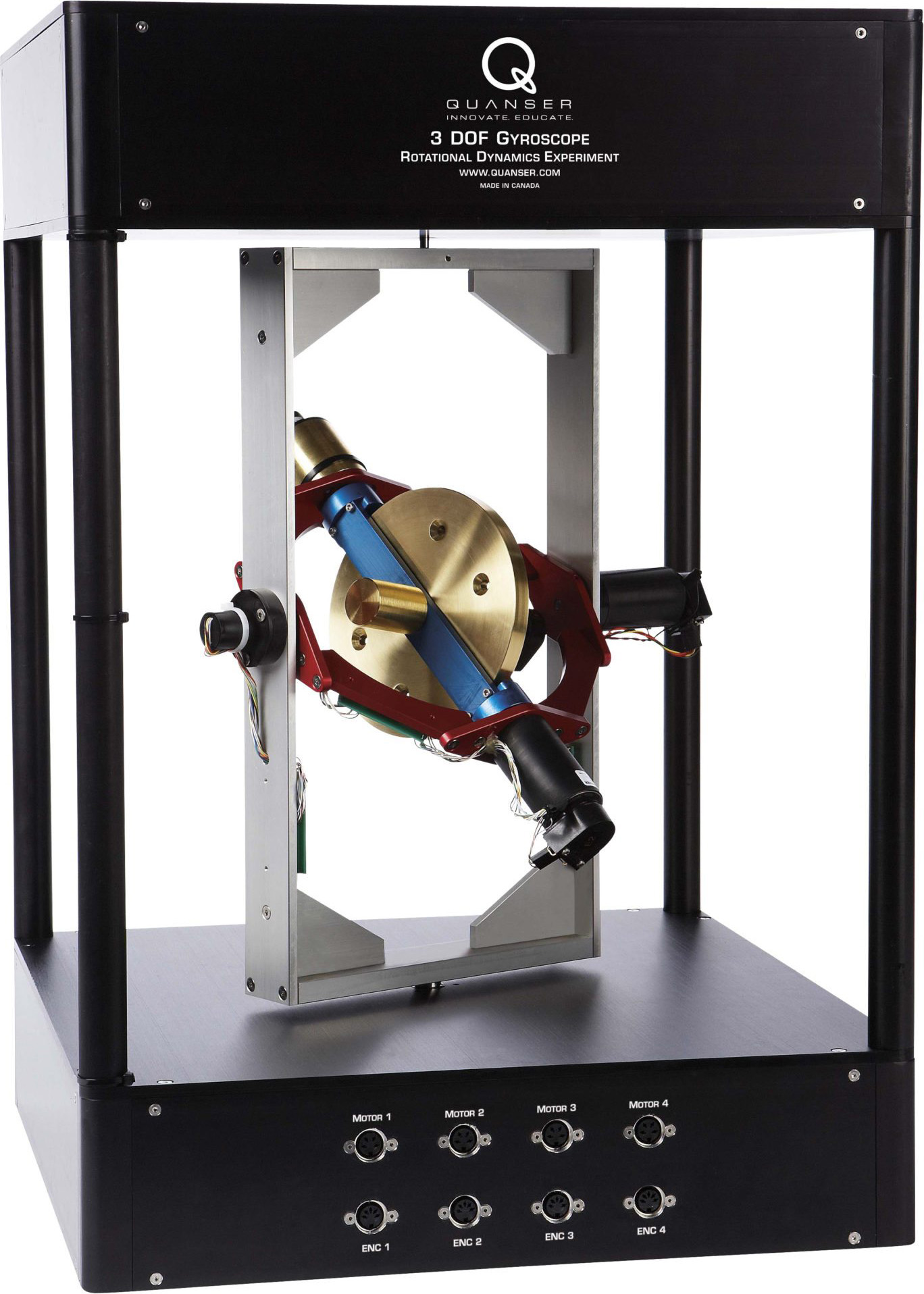} 
\caption{CMG system} \label{fig:gyro}
\end{subfigure}
\caption{Considered benchmark systems, whose  
simulation model is used for identification.} \vspace{-3mm}
\label{fig:benchmarks}
\end{figure}

In this section we demonstrate the estimation performance and computational efficiency of the proposed identification scheme on two benchmark systems, depicted in Figure~\ref{fig:benchmarks}, utilizing simulation data. Estimation of LTI, LPV and general NL models will be accomplished, showcasing the flexibility of the overall system identification toolchain.  

\subsection{Implementation and computational aspects}
The proposed identification approach is implemented with the {\tt jax-sysid} package~\cite{Bem25}
in Python 3.11 and run on an MacBook Pro with Apple M4 Max processor equipped with 16 CPU cores. Given an initial guess of the unknowns $(\theta,\hat{w}_0)$, we report the CPU time spent on solving Problem~\eqref{eq:optim-prob} on a single core.

In all the examples, unless specified differently, we set: the regularization
term \eqref{eq:regularization} to $R(\theta,\hat{w}_0) = \frac{\rho_\theta}{2}\|\theta \|_2^2
+ \frac{\rho_\mathrm{w}}{2}\|\hat{w}_0 \|_2^2$
with $\rho_\theta=2\cdot10^{-4}$ and $\rho_\mathrm{w}=2\cdot10^{-8}$; if Adam is used
for warm start, we run $10^3$ iterations with parameters $\eta=10^{-3}$, $\beta_1=0.9$, $\beta_2=0.99$, $\epsilon=10^{-8}$; L-BFGS-B optimization is run for at most $10^4$ iterations;  
the results shown in the tables are obtained by repeating the identification procedure 100 times from different random initial guesses, selecting the model providing the highest \emph{best fit rate} (BFR)\footnote{For prediction, $\mathrm{BFR}=\mathrm{max}\left\{1-\sqrt{\frac{\sum_{k=0}^{N-1} \| y_k - y_k^\mathrm{pred} \|_2^2 }{\sum_{k=0}^{N-1} \| y_k - \bar{y}  \|_2^2}},0\right\}$ with $\bar{y}$ being the sample mean of the measured $y$, while for simulation, $y_k^\mathrm{pred}$ is substituted with $y_k^\mathrm{sim}$, resulting when \eqref{eq:1-step-prediction} is run with $\tilde f_\mathrm{z}$ and $\tilde g_\mathrm{z}$ set to $0$.} of the predicted and simulated model response on test data in terms of \eqref{eq:1-step-prediction}; unless specified differently, the coefficients of the deterministic model~\eqref{eq:model-x} are not used as a warm start when training the combined model.

For each example, we identify the combination of plant~\eqref{eq:model-x} and (inverse) noise model~\eqref{eq:model-z-inv} in the form \eqref{eq:1-step-prediction} by solving~\eqref{eq:optim-prob}, denoted as "{\combined}"
in the reported tables. For comparison, we also train a purely deterministic plant model (\eqref{eq:model-x} only)  by setting $\tilde f_\mathrm{z}$ and $\tilde g_\mathrm{z}$ to $0$ in~\eqref{eq:optim-prob}, denoted as "\plantonly" in the tables.
Note that \eqref{eq:1-step-prediction} corresponds to the simulation error when only a plant model is trained, realizing an OE noise-based model structure.  After identifying a model, the initial state is reconstructed from
the dataset as described in~\cite[Section III.C]{Bem25}
to compute output predictions. The tables also show the maximum achievable performance, denoted as "{\tt true}"\footnote{Results under "{\tt true}" mean that the original representation $(G_\mathrm{o},H_\mathrm{o})$ of the data-generating system has been used as a model to calculate (i) the simulation performance in terms of subtracting the noise free response $y_{\mathrm{o},k}$ of the process model $G_\mathrm{o}$ from the measured output $y$, which gives back the true $v_k$ as the error in the BFR calculation; (ii) the prediction performance, which is simply $e_k=H_\mathrm{o}^{-1} v_k$, i.e., the original noise sequence is used as the prediction error in the BFR calculations. }.
The CPU time reported in the tables is referred to the training of the best model. 

\subsection{Unbalanced disk}
\label{benchmark:unbdisk}
\subsubsection{Data-generating system} The unbalanced disk system, depicted in Fig.~\ref{fig:unbd}, consists of a vertically mounted disc with a lumped mass $m$ at distance $l$ from the center, actuated by an input voltage $u$. Due to the gravitational force exerted on the lumped mass, the system exhibits nonlinear dynamics, which can be represented by the following differential equation:
\begin{equation}
	\ddot{\alpha} = -\frac{1}{\tau} \dot{\alpha} + \frac{K_\mathrm{m}}{\tau} u - \frac{m g l}{J} \sin(\alpha),
	\label{eqn:unb_disc_nonlinear}
\end{equation}
where $\alpha$ is the disc angular position, $\tau$ is the lumped back EMF constant, $K_\mathrm{m}$ is the motor constant, $J$ is the complete disc inertia, and $g$ is the acceleration due to gravity. These physical parameters are given in \cite{Toth21Sysidac}. To generate discrete-time data, under \emph{zero-order-hold} (ZOH) actuation, \eqref{eqn:unb_disc_nonlinear} is discretized via the Euler approach with sampling time $T_\mathrm{s}=0.01$s. The measured output is the angular position $\alpha$ corrupted by an additive noise process $v_k$, that will vary in the different scenarios that we will consider.  

\subsubsection{LTI case}
To test our approach in the ideal case of data generated by a linear system, we linearize~\eqref{eqn:unb_disc_nonlinear} 
around $(\alpha,\dot{\alpha})=(0,0)$ and discretize the resulting linear dynamics, obtaining the plant dynamics  \vspace{-1mm}
\begin{subequations}
\label{eq:LTI-unbalanced}
\begin{align}
x_{k+1}& =A_\mathrm{x} {x}_k+B_\mathrm{x} u_k, \\
y_{\mathrm{o},k} & = C_\mathrm{x}   x_k,
\end{align}
with \vspace{-2mm}
\begin{equation}
\begin{aligned}
	A_{\mathrm{x},0} &=\begin{bmatrix} 1-\frac{T_\mathrm{s}}{\tau} &  0 & 1\end{bmatrix}\!, \ A_{\mathrm{x},1}=\begin{bmatrix} 0 &  -\frac{Mg\ell T_\mathrm{s}}{J}\\0 & 0\end{bmatrix}\!,\\
	A_\mathrm{x}&=A_0 +A_1, \
	B_\mathrm{x}=\begin{bmatrix} \frac{T_\mathrm{s}K_\mathrm{m}}{\tau}\\0\end{bmatrix}\!,\
	C_\mathrm{x}=\begin{bmatrix} 0 & 1\end{bmatrix}\!,
\label{eq:unbalanced-matrices}
\end{aligned}
\end{equation}
completed with the BJ noise dynamics
\begin{align}
z_{k+1}&=(A_{\mathrm{z},0}+A_{\mathrm{z},1})z_{k}+(B_{\mathrm{z},0}+B_{\mathrm{z},1})e_k, \\
v_k &=z_k+e_k, \\
y_k &= y_{\mathrm{o},k}+ v_k, 
\end{align}
\end{subequations}
where $A_{\mathrm{z},0}=\frac{277}{300}$, $A_{\mathrm{z},1}=\frac{2}{30}$, $B_{\mathrm{z},0}=\frac{41}{150}$, $B_{\mathrm{z},1}=\frac{-11}{60}$.
Using \eqref{eq:LTI-unbalanced}, independent training and test datasets, $N=2000$ sample long each, are generated by using  $u_k\sim{\mathcal N}(0,1)$ and $e_k\sim{\mathcal N}(0,3.75\cdot 10^{-3})$, giving a \emph{signal-to-noise ratio} (SNR) of about 10~dB on both sets. 
Using our proposed approach, we run L-BFGS-B, without the need of Adam warm starting, and identify an LTI plant model with $\hat{n}_\mathrm{x}=2$ states and $\hat{n}_\mathrm{z}=1$ noise-model state, corresponding to the minimal state dimensions required to describe the model dynamics in terms of \eqref{eq:LTI-unbalanced}. For comparison purposes,  the subspace method {\tt n4sid} and the gradient-based PEM-SS approach {\tt ssest} from the System Identification Toolbox in {\textsc{Matlab}}~\cite{Lju01} are also applied on the datasets with default settings.

The results are summarized in Table~\ref{tab:LTI-unbalanced}, where the column ``type'' refers to whether the BFR calculations are based on the predicted or simulated response of the model. The labels {\tt (s)} and {\tt (p)} for {\tt ssest} indicate if the estimation is performed with simulation or prediction focus, respectively. Note that to capture only the process dynamics, $\hat{n}_\mathrm{x}=2$ is sufficient, while to capture also the noise dynamics $\hat{n}_\mathrm{x}=3$ is required in an innovation form used by the subspace and PEM-SS methods in \textsc{Matlab}, hence estimation is performed with the required optimal state dimensions for these methods, showing the obtained results in Table~\ref{tab:LTI-unbalanced}. Based on our results, \texttt{n4sid} with its default options fails to find a good model in this case, while \texttt{ssest} achieves worse results in term of the estimation of the process part than our proposed approach, with comparable 
computation time.

\begin{table}[h!]
\begin{center}
\caption{Identification results for the linearized unbalanced disk system~\eqref{eq:LTI-unbalanced} with BJ noise dynamics. }
\label{tab:LTI-unbalanced}
\begin{tabular}{@{}l|c|c|r|r|c|r@{}}
  & $\hat{n}_\mathrm{x}$ & $\hat{n}_\mathrm{z}$ & BFR train. & BFR test & type & time\\
  \hline
  \multirow{2}{*}{\tt true} & \multirow{2}{*}{2} & \multirow{2}{*}{1} & 72.12\% & 68.13\% & sim & - \\
  & & &  77.66\% & 72.85\% & pred & - \\ 
\hline
{\plantonlyt} &
2 & 0 &  72.03\% &  68.08\% & sim &  0.13\sps s\\
\hline
\multirow{2}{*}{\combinedt} & 2 & 1 &  72.02\% &  68.08\% & sim & \multirow{2}{*}{0.34\sps s}\\
& 2 & 1 &  77.67\% &  72.84\% & pred\\
\hline
{\tt n4sid} {\tt (s)} & 2 & - & 0.33\% & 0.56\% & sim &  0.15\sps s\\
\multirow{2}{*}{{\tt n4sid}  {\tt (p)}} & \multirow{2}{*}{3} & \multirow{2}{*}{-} & 61.21\% & 54.49\% & sim &  \multirow{2}{*}{0.11\sps s}\\
& & &  75.63\% & 70.46\% & pred \\
{\tt ssest} {\tt (s)} & 2 & - & 1.33\% & 1.37\% & sim &  0.47\sps s\\
\multirow{2}{*}{{\tt ssest} {\tt (p)}} & \multirow{2}{*}{3} & \multirow{2}{*}{-} & 64.23\% & 58.02\% & sim &  \multirow{2}{*}{0.20\sps s}\\
& & &  76.31\% & 71.43\% & pred \\
\hline
\end{tabular}
\end{center}
\end{table}

\subsubsection{Externally scheduled LPV case} 
As discussed in \cite{Toth21Sysidac}, we convert the nonlinear dynamics of the
unbalanced disk system to  an LPV form that we discretize, giving the LPV plant dynamics combined with a scheduling-dependent noise process in the following form:
\begin{subequations}
\label{eq:LPV-unbalanced}
\begin{align}
x_{k+1}& =\overbrace{(A_{\mathrm{x},0}+A_{\mathrm{x},1}p_k)}^{A_\mathrm{x}(p_k)} {x}_k+B_\mathrm{x} u_k,\\
z_{k+1}&=(A_{\mathrm{z},0}+A_{\mathrm{z},1}p_k)z_{k}+(B_{\mathrm{z},0}+B_{\mathrm{z},1}p_k)e_k,\\
y_{k} & = C_\mathrm{x}   x_k+ \underbrace{z_k+e_k}_{v_k},
\end{align}
\end{subequations}  
with $A_{\mathrm{x},0},\ldots,C_\mathrm{x}$ and $A_{\mathrm{z},0},\ldots,B_{\mathrm{z},1}$ defined as in~\eqref{eq:unbalanced-matrices}.
We assume that the scheduling signal $p_k$ is known and exogenously generated,
$p_k=(1-p_{\rm mag})+p_{\rm mag}n_k$, where $p_{\rm mag}=0.25$ and $n_k$ is a random binary signal with bandwidth $\frac{0.05}{2T_s}$~Hz, while the input $u$ and the generating noise $e$ are independent white noise sequences with $u_k\sim{\mathcal N}(0,1)$ and $e_k\sim{\mathcal N}(0,3.75\cdot 10^{-3})$ corresponding to an SNR of 10 dB. 
Under these conditions, a training and test dataset are generated independently, with $N=2000$ samples each. Using the proposed identification approach, we identified models with $\hat{n}_\mathrm{x}=2$ and $\hat{n}_\mathrm{z}=1$, while all matrix functions in \eqref{eq:LPV-x} and \eqref{eq:LPV-z-inv} are assumed to depend affinely on $p_k$. 
For comparison, we used the {\textsc{LPVcore}} Toolbox for \textsc{\textsc{Matlab}}~\cite{Toth21Sysidac}
to identify an LPV state-space model of the same order 
using the subspace method {\tt lpvidss} with a past window of $6$ and the gradient-based PEM-SS method {\tt lpvssest}.

Table~\ref{tab:LPV-unbalanced} shows the results obtained by running the L-BFGS-B optimization
without Adam warm starting. As it can be seen, the estimation performance of the newly proposed approach is slightly better than gradient-based PEM in \textsc{LPVcore} and, as expected, achieves better results than the subspace method, due the direct minimization of the $\ell_2$-loss on the prediction error. In terms of computation time, the subspace approach is faster due to the short past window taken, while the newly proposed method, implemented in Python with JAX, for a given initial guess of the parameters,
converges faster than {\tt lpvssest} in \textsc{\textsc{Matlab}}. 

\begin{table}[h!]
\begin{center}
\caption{Identification results for the LPV form of the unbalanced disk system~\eqref{eq:LPV-unbalanced} with LPV-BJ noise dynamics and external scheduling.}
\label{tab:LPV-unbalanced}
\begin{tabular}{@{}l|c|c|r|r|r|c|r@{}}
  & $\hat{n}_\mathrm{x}$ & $\hat{n}_\mathrm{z}$ & sched. & BFR train. & BFR test & type & time\\
   \hline
  \multirow{2}{*}{\tt true} & \multirow{2}{*}{2} & \multirow{2}{*}{1} &   \multirow{2}{*}{ ext.}& 70.69\% & 72.61\% & sim & - \\
  & & & &  73.63\% & 74.36\% & pred & - \\ 
\hline
\plantonlyt & 2 & 0 & ext. & 70.84\% &   71.94\% & sim &  3.81\sps s\\
\hline \multirow{2}{*}{\combinedt}
& \multirow{2}{*}{2}  & \multirow{2}{*}{1}  & ext.&70.82\% &  71.81\% & sim & \multirow{2}{*}{5.15\sps s}\\
&  &  & ext. &73.62\% & 74.19\% & pred \\
\hline
\multirow{2}{*}{\tt lpvidss} & \multirow{2}{*}{3} & \multirow{2}{*}{-} & \multirow{2}{*}{ext.} & 1.33\% & 1.37\% & sim &  \multirow{2}{*}{1.99\sps s}\\
& & & & 67.45\% & 66.81\% & pred \\ \hline
\multirow{2}{*}{\tt lpvssest} & \multirow{2}{*}{3} & \multirow{2}{*}{-} & \multirow{2}{*}{ext.} & 69.80\%  & 69.48\% & sim & \multirow{2}{*}{21.78\sps s}\\
& & & & 73.37\% & 72.39\% & pred  \\
\hline
\end{tabular}
\end{center}
\end{table}

\subsubsection{Self-scheduled LPV case} 
We now identify a self-scheduled LPV model with $\hat{n}_\mathrm{x}=2$ and $\hat{n}_\mathrm{z}=1$, while all matrix functions in \eqref{eq:LPV-x} and \eqref{eq:LPV-z-inv} are assumed to depend on $p_k$ in terms of affine dependence. To generate data, \eqref{eq:LPV-unbalanced} has been simulated with $p_k=\mathrm{sinc}(x_2(k))= \frac{\sin(\alpha(k))}{\alpha(k)}$, corresponding to the original nonlinear dynamics of the unbalanced disk system. Under the same data generation with $N=2000$, $u_k\sim{\mathcal N}(0,1)$, $e_k\sim{\mathcal N}(0,3.75\cdot 10^{-3})$, the obtained estimation and test data sets had an SNR of 21 dB.
For the scheduling map $\psi$, a feedforward ANN with two hidden layers is used, with 6 nodes in each and \emph{sigmoid} (first layer) and \emph{swish} (second layer) activation functions.
In total, the plant model has 83 parameters to learn,
while the combined plant and noise model has 86 (LTI noise model)
and 89 (LPV noise model) parameters.
The results achieved using 1000 Adam iterations to warm start the L-BFGS optimization,
are summarized in Table~\ref{tab:qLPV-unbalanced}. As we can see, an LPV noise model 
provides better results than an LTI noise model.

\begin{figure}[h!]
\centering
\begin{subfigure}[b]{\textwidth}
    \centering
    \includegraphics[scale=0.60, trim=25 40 40 80, clip]{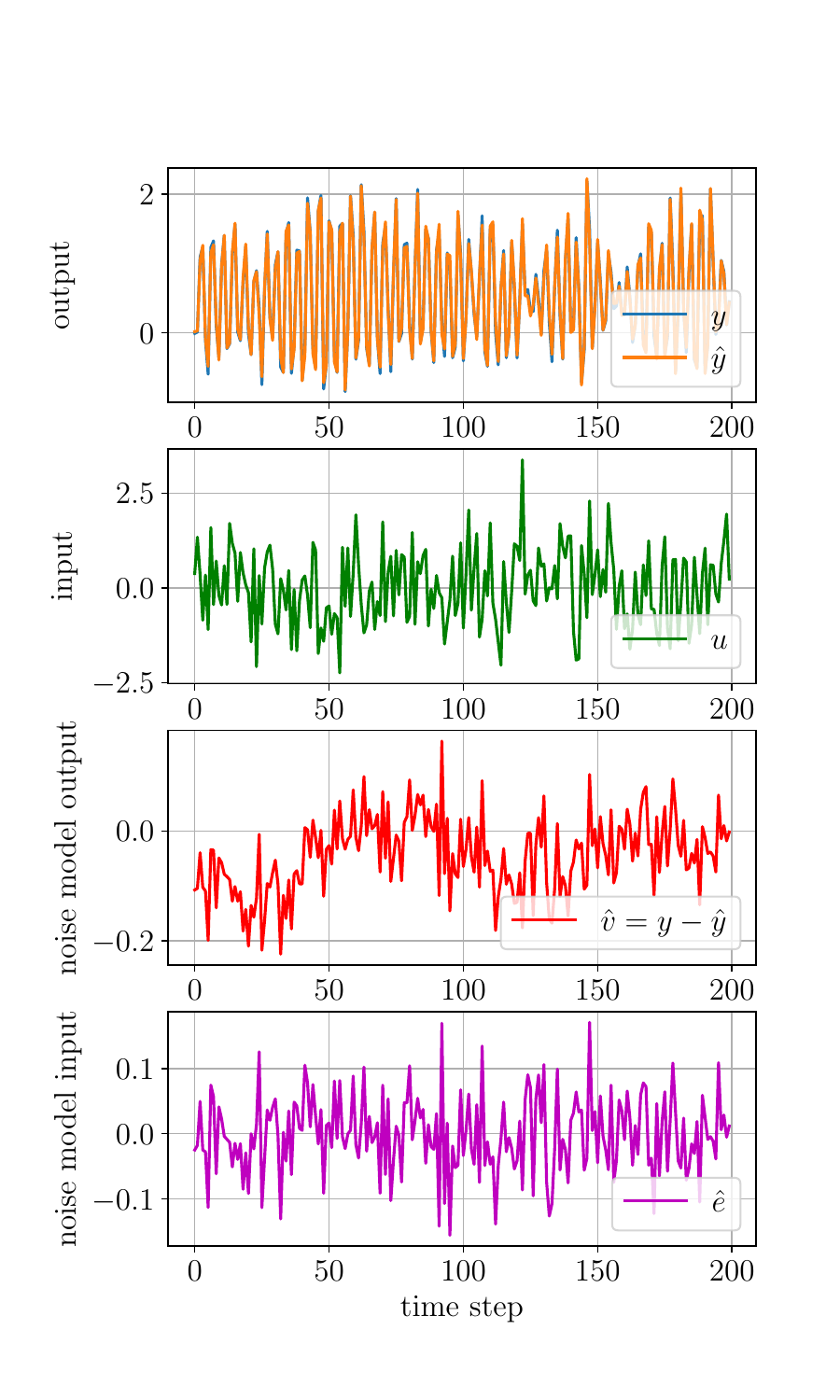}
    \caption{First 200 samples of one-step-ahead predictions} \vspace{-3mm}
    \label{fig:qLPV-unbalanced-time}
\end{subfigure}
\\\vspace{3em}
\begin{subfigure}[b]{\textwidth}
    \centering
    \includegraphics[scale=0.35, trim=20 40 60 20,clip]{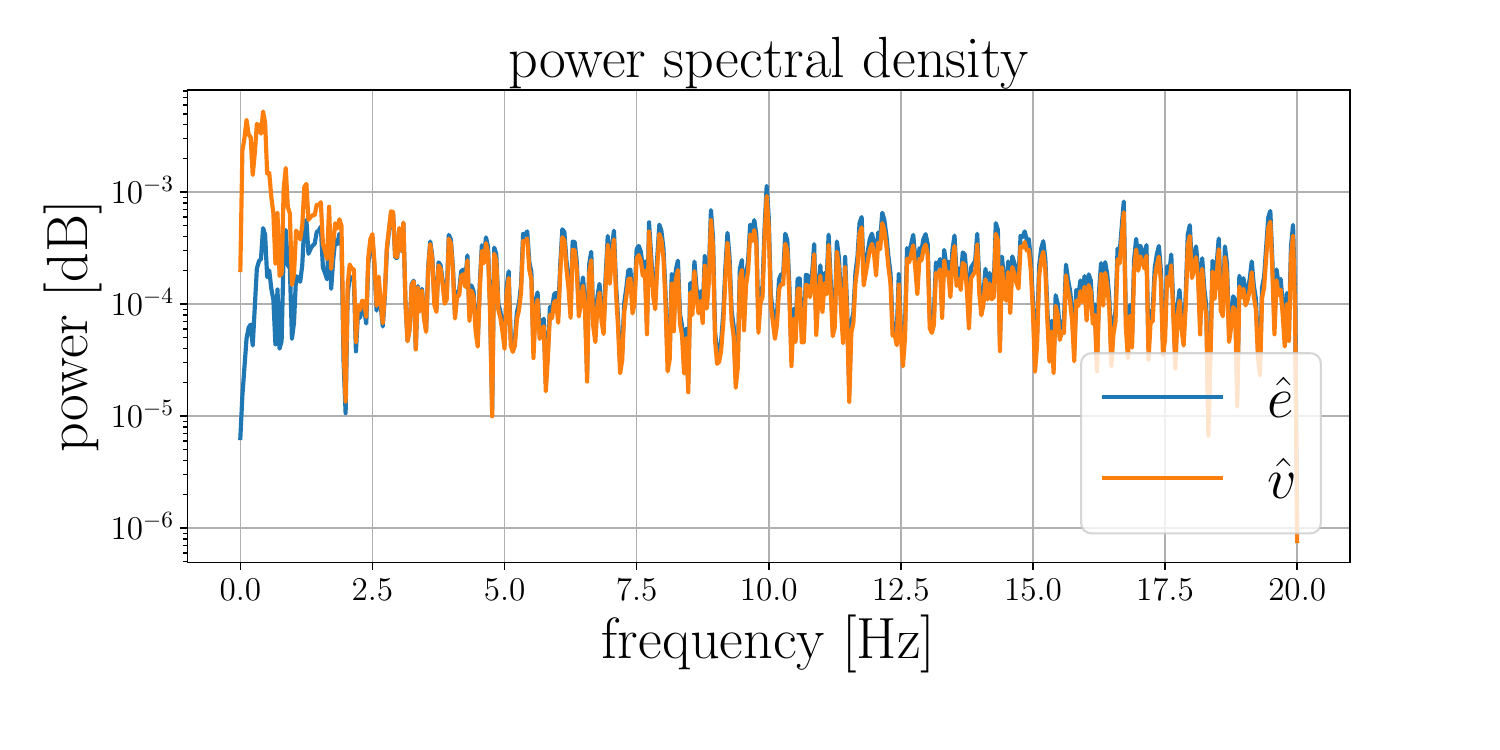}
    \caption{Power spectral density of the reconstructed noise $\hat e_k$ and $\hat v_k$}
    \label{fig:qLPV-unbalanced-freq}
\end{subfigure}
\caption{Prediction results of the self-scheduled LPV model with LTI noise model identified for the unbalanced disk.}
\label{fig:qLPV-gyroscope}
\end{figure}

\begin{table}[h!]
\setlength{\tabcolsep}{5pt}
\begin{center}
\caption{Identification results for the unbalanced disk system with input noise and using self-scheduled LPV models.}
\label{tab:qLPV-unbalanced}
\begin{tabular}{@{}l|c|c|c|r|r|c|c|c@{}}
  & $\hat{n}_\mathrm{x}$ & $\hat{n}_\mathrm{z}$ & sched. & BFR train. & BFR test & type & noise  & time\\
\hline
  \multirow{2}{*}{\tt true} & \multirow{2}{*}{2} & \multirow{2}{*}{1} &   \multirow{2}{*}{self}& 89.32\% & 90.11\% & sim & LPV & - \\
  & & & &  91.24\% & 91.86\% & pred & LPV & - \\ 
     \hline
\plantonlyt &
2 & 0 & self & 87.73\% &  86.45\% & sim & -- & {10.06\sps s}\\
\hline
  \multirow{2}{*}{\combinedt} &  \multirow{2}{*}{2} &  \multirow{2}{*}{1} &  \multirow{2}{*}{self} & 85.19\% &  86.19\%  & sim& LTI & \multirow{2}{*}{12.91\sps s}\\
&  &  & & 90.92\% &   91.46\% & pred& LTI\\
\hline
\multirow{2}{*}{\combinedt} &\multirow{2}{*}{2} & \multirow{2}{*}{1} & \multirow{2}{*}{self} & 85.60\% &  86.56\% & sim & LPV & \multirow{2}{*}{18.82\sps s}\\
& & &  & 90.96\% &  91.51\% & pred & LPV\\
\hline
\end{tabular}
\end{center}
\end{table}
For the case of identification of a self-scheduled LPV model (\combined), the sample variance of $\hat v_k$ and $\hat e_k$ on test data are, respectively,
$\Var\{\hat v_k\} = 5.5981\times 10^{-3}$, $\Var\{\hat e_k\} = 4.3937\times 10^{-3}$.
Figure~\ref{fig:qLPV-unbalanced-time} shows the one-step-ahead predictions obtained by 
the \combined\ model with LTI noise model for the first 200 samples of the test data. 
It is apparent how the noise model has tried to capture the acceleration noise.
The plot also shows the evolution $\hat v_k$ and $\hat e_k$, while their power spectral density is shown in Figure~\ref{fig:qLPV-unbalanced-freq}, where it can be observed
that the reconstructed noise $\hat e_k$ has less low-frequency content than the estimate $\hat v_k$
of the output error $v_k=y_k-y_{\mathrm{o},k}$. 

\subsubsection{Nonlinear case}
We identify a nonlinear state-space plant model~\eqref{eq:model-x} in which $f_\mathrm{x}$ is a two-layer
ANN with 15 and 10 neurons and \emph{swish} activation function, while
$f_\mathrm{y}$ is a linear function. The noise model is taken to be LTI with $\hat{n}_\mathrm{z}=1$ state.
The total number of training parameters is 244 (plant model only) and 247 (plant + noise model). The  results are summarized in Table~\ref{tab:NL-unbalanced}. We can see that the proposed approach can be easily used to scale up the model estimation to the NL dynamics, and can reliably capture the underlying behavior of the system. 
Note that the simulation performance on training data of \plantonly~ is slightly better (+0.01\%) than {\tt true}; this is possibly due to the extra degree of freedom the training algorithm has in optimizing also the initial state $\hat x_0$. In fact, if we increase the weight on $\hat x_0$ to $2\cdot10^{-3}$, the BFR on training data becomes 89.32\%.

\begin{table}[h!]
\begin{center}
\caption{Identification results for the altered unbalanced disk system with input noise and using NL models.}
\label{tab:NL-unbalanced}
\begin{tabular}{@{}l|c|c|r|r|c|r@{}}
  & $\hat{n}_\mathrm{x}$ & $\hat{n}_\mathrm{z}$ & BFR train. & BFR test & type & time\\
\hline
  \multirow{2}{*}{\tt true} & \multirow{2}{*}{2} & \multirow{2}{*}{1} & 89.32\% & 90.11\% & sim  & - \\
  & & &  91.24\% & 91.86\% & pred  & - \\ 
     \hline
\plantonlyt & 2 & 0 &  89.33\% &  89.89\%& sim &  10.24\sps s\\
\hline
  \multirow{2}{*}{\combinedt} & \multirow{2}{*}{2} & \multirow{2}{*}{1} &  89.23\%& 89.83\%& sim & \multirow{2}{*}{11.31\sps s}\\
&  &  &  91.05\% &  91.39\%& pred\\
\hline
\end{tabular}
\end{center}
\end{table}

\subsection{Control moment gyroscope}
\label{benchmark:cmg}
We consider a high-fidelity NL simulator of the 3 degrees-of-freedom \emph{Control Moment Gyroscope} (CMG) depicted in Fig.~\ref{fig:gyro}. 

\subsubsection{Data-generating system} \label{benchmark:cmg:datagen} The system consists of a flywheel (with angular position $q_1$), mounted inside an inner blue gimbal (with angular position $q_2$), which in turn is mounted inside an outer red gimbal (with angular position $q_3$). The entire structure is supported by a rectangular silver frame (with angular position $q_4$) that can rotate around its vertical axis of symmetry. The gimbals can rotate freely. The setup is equipped with four DC motors and encoders, actuating and measuring the position of the flywheel and the gimbals, respectively. The nonlinear system is highly-complex due to the involved rotational dynamics. 

In this study, we use the high-fidelity CMG simulation model~\cite{bloemers2019equations} as the data-generating system, which is implemented in C and interfaced to \textsc{Matlab}. We consider a scenario in which the red gimbal is locked ($q_3\equiv 0$), the blue gimbal ($q_2$) is actuated with motor current $i_2$ (treated as the manipulated input), while $\dot{q}_4$ is the measured output, corrupted by a white noise signal $e_k\sim \mathcal{N}(0,\sigma_\mathrm{e}^2) $ with variance $\sigma_\mathrm{e}^2=2.2 \cdot 10^{-5}$. The angular velocity $\dot{q}_1$ of the flywheel is controlled independently to follow a random multi-level-reference signal with amplitude in $[30,50]$ rad/s and a dwell-time between 4 and 8 seconds. 

The CMG is simulated in continuous time with fixed-step Runge-Kutta 4 under ZOH, sampling the results every $T_\mathrm{s}=0.01$~s to obtain the training and test data sets. By analyzing the first-principle dynamics of the system (see, e.g., \cite{bloemers2019equations}), it can be noted that, due to the independent control of the flywheel and the locked status of $q_3$, the angular velocity $\dot{q}_1$ can be regarded as an external input to a model describing the motion of $q_2$ and $q_4$. Hence, in the considered scenarios, $\dot{q}_1$ is treated either as an additional input signal (LTI, self-scheduled LPV, and NL case) or a scheduling variable (externally-scheduled LPV case). 

\subsubsection{LTI case} 
\label{sec:CMG:LTI}
Using the simulator with the configuration described in Section \ref{benchmark:cmg:datagen}, a training dataset of size $N=10^4$ and test dataset of size $N=3\cdot 10^4$ were generated by exciting $i_2$ with the experiment design detailed in \cite[Sec.~IV.B]{Toth22CDCb}, giving an SNR of 35 dB on both data sets.  
We identify LTI models of the nonlinear dynamics inputs $u_1=i_2$, $u_2=\dot q_1$.

Using model-structure selection, the plant model order $\hat{n}_\mathrm{x}=8$ was established as a good choice while  $\hat{n}_\mathrm{z}=2$ was taken to test overfitting as, in this case, the data-generating system has an OE noise structure, meaning that the estimation can not benefit from a noise model. For comparison purposes, {\tt n4sid} with focus on simulation and {\tt ssest} with focus on simulation {\tt (s)} and prediction {\tt (p)} were also executed on the data set using $\hat{n}_\mathrm{x}=8$ (corresponding to \plantonly) and $\hat{n}_\mathrm{x}=10$ (corresponding to \combined) as model order. The results, obtained by running L-BFGS without Adam warm start and $\rho_\theta=2\cdot10^{-9}$,
are shown in Table~\ref{tab:LTI-gyroscope}.

As can be seen from the table, the proposed approach results in a better fit on the training data than {\tt n4sid} or {\tt ssest} when only a plant model is estimated, as well as on test data, except for {\tt ssest} {\tt (p)}. When both a plant and a noise model are estimated, the best fit on test data is achieved by the proposed method, without excessively sacrificing the fit on training data.

\begin{table}[h!]
\begin{center}
\caption{Identification results for the CMG system (Section \ref{benchmark:cmg})  with LTI models.}
\label{tab:LTI-gyroscope}
\begin{tabular}{@{}l|c|c|r|r|c|r@{}}
  & $\hat{n}_\mathrm{x}$ & $\hat{n}_\mathrm{z}$ & BFR train. & BFR test & type & time\\
  \hline
  {\tt true} &  5  & -  & 98.16\% &  98.19\% & sim & - \\
\hline
\plantonlyt & 8 & 0 &  35.99\%&  25.28\% & sim &  3.50\sps s\\
\hline
\multirow{2}{*}{\combinedt}
& 8 & 2 &  34.46\% &  29.91\% & sim & \multirow{2}{*}{8.73\sps s}\\
& 8 & 2 &  97.17\% &  97.12\% & pred\\
\hline
{\tt n4sid} & 8 & 0 & 29.72\% & 20.98\% & sim &  0.85\sps s\\
{\tt n4sid} & 10 & 0 & 34.76\% & 22.45\% & sim &  1.05\sps s\\
{\tt ssest} {\tt (s)} & 8 & 0 & 35.48\% & 24.70\% & sim &  2.05\sps s\\
{\tt ssest} {\tt (p)} & 10 & 0 & 33.51\% & 26.88\% & sim &  3.06\sps s\\
\hline
\end{tabular}
\end{center}
\end{table}

\subsubsection{LPV case} 
\label{sec:CMG:LPV}
We use the same simulation setting to identify LPV models.
The plant model has $\hat{n}_\mathrm{x}=5$ states and affine dependence on $p$, while the noise model is taken LTI with $\hat{n}_\mathrm{z}=2$ states.  
To make a fair comparison with results previously obtained on this dataset in~\cite{Toth22CDCb}, for the self-scheduled LPV case we choose $n_\mathrm{p}=3$ and parameterize the scheduling map $\psi$ as a feedforward ANN with two hidden layers of 5 neurons each with a  
\emph{swish} activation function. The plant and the combined plant $+$ noise models have 281 and 289 parameters, respectively. For the externally-scheduled case, 
the scheduling variable $p_k$ is taken from measured data by setting $p = \psi(x,u) = [\dot{q}_1\ \sin(q_2)\ \cos(q_2)]^\top$ and the plant and the combined plant $+$ noise models have 185 and 193 parameters, respectively. 

For initialization of the proposed identification scheme, we use a bootstrapping process by first training an LTI plant model, use it to initialize the training of an LPV plant model (\plantonly), which is in turn used to initialize the training of the full model including a noise model (\combined). 
Note again that, due the OE noise structure of the data-generating system, the noise model is taken only to test overfitting by the proposed estimation scheme. The results obtained by training the models, with $\rho_\mathrm{w}=2\cdot10^{-4}$ 
and $\rho_\theta=2\cdot10^{-8}$ (externally-scheduled case), and $\rho_\theta=2\cdot10^{-3}$ (self-scheduled case),
are summarized in Table~\ref{tab:LPV-gyroscope}. For comparison purposes, the table also shows the results obtained in \cite{Toth22CDCb} using the input-output model structure based PEM approach {\tt lpvoe} with a $5^\mathrm{th}$ order LPV-OE input-output model ($n_\mathrm{a}=n_\mathrm{b}=5$) and the gradient-based PEM-SS method {\tt lpvssest} with $\hat{n}_\mathrm{x}=5$ and no noise model estimation implemented in \textsc{LPVcore}. 
Next to this, the LPV variant of the SUBNET method has been also applied to identify a self-scheduled LPV model with $\hat{n}_\mathrm{x}=5$ and $n_\mathrm{p}=3$ using a scheduling map parameterization with 2 hidden layers, each with 64 neurons, $\tanh$ activation functions, and a linear bypass. 

As we can see from Table~\ref{tab:LPV-gyroscope}, the proposed approach with only process model estimation achieves better results than the state-of-the-art SUBNET approach, both for the self and externally scheduled cases. Actually, it comes nearly to the performance of the true model (theoretical maximum) in the externally scheduled case and accomplishes the estimation in a fraction ($\approx$\sps 0.2\sps\%) of the required training time of SUBNET, which is a remarkable achievement. Note that no other approach is currently capable of combined estimation of the scheduling map and the LPV model than the proposed method and the SUBNET approach. When compared to traditional LPV approaches working with a user defined scheduling map such as {\tt lpvoe} and {\tt lpvssest}, it is obvious that co-estimation of the scheduling map and the plant model gives clear advantages for both the proposed method and SUBNET.

When a noise model is added, in the self-scheduled case the estimated model results in worse simulation performance compared to the \plantonly ~case. However, in case of external scheduling, the performance is close to that of the \plantonly ~case and the noise model is consistently estimated to be approximately zero.    

\begin{table}[h!]
\begin{center}
\caption{Identification results for the CMG system (Section \ref{benchmark:cmg})  with LPV models.}
\label{tab:LPV-gyroscope}
\begin{tabular}{@{}l|c|c|c|r|r|c|c@{}}
  & $\hat{n}_\mathrm{x}$ & $\hat{n}_\mathrm{z}$ & sched. & BFR train. & BFR test & type & time\\
   \hline
  {\tt true} &   5 & -  & - &  98.16\% &  98.19\% & sim & - \\
\hline
\multirow{2}{*}{\plantonlyt} &\multirow{2}{*}{5} & \multirow{2}{*}{0} & self &  97.61\% &  96.50\% & \multirow{2}{*}{sim} &  47.54\sps s\\
                            &    &  & ext. &  98.14\% &  98.13\% &  &  22.65\sps s\\
\hline
\multirow{2}{*}{\combinedt} &\multirow{2}{*}{5} & \multirow{2}{*}{2} & \multirow{2}{*}{self} &81.96\% & 83.78\% & sim & \multirow{2}{*}{67.19\sps s}\\
& &  &   & 97.56\%  & 97.64\% & pred\\
\hline
\multirow{2}{*}{\combinedt} &\multirow{2}{*}{5} & \multirow{2}{*}{2} & \multirow{2}{*}{ext.} &98.12\% & 98.11\% & sim & \multirow{2}{*}{31.61\sps s}\\
& &  &   & 98.11\%  & 98.10\% & pred\\
\hline
{\tt lpvoe} & 5 & 0 & ext. & 93.32\% & 89.96\% & sim & 11.12\sps s \\ \hline
\multirow{2}{*}{\tt lpvssest} & \multirow{2}{*}{5} & \multirow{2}{*}{0} & \multirow{2}{*}{ext.} & 93.62\% & 90.84\% & sim &  \multirow{2}{*}{348\sps s}\\ 
& &  &   & 94.91\%  & {93.22}\% & pred \\
\hline
\multirow{2}{*}{SUBNET}& \multirow{2}{*}{5} & \multirow{2}{*}{0} & self. &  {97.28}\%& {96.40}\% & \multirow{2}{*}{sim} & $\approx$\sps 10\sps h  \\
 &  &  & ext. &  {97.25}\% & {97.01}\% &  & $\approx$\sps 10\sps h \\\hline
\end{tabular}
\end{center}
\end{table}

\subsubsection{Nonlinear case}
\label{sec:NL-gyroscope}
Finally, we identify a fully nonlinear model using the same simulator and configuration for the CMG system described above.
We let $f_\mathrm{x}$ be the sum of a linear function $Ax_k+Bu_k$ and
a two-layer FNN with 10 neurons each and $f_\mathrm{y}$ the sum
of a linear function $Cx_k$ and a two-layer FNN with 5 neurons each, both with 
\emph{swish} activation functions, with $\hat{n}_\mathrm{x}=5$ states, and set $\rho_\theta=2\cdot10^{-3}$, $\rho_\mathrm{w}=2\cdot10^{-4}$. The noise model is LTI with $\hat{n}_\mathrm{z}=2$ states. Again, the same bootstrapping process is used for initialization as discussed in Section \ref{sec:CMG:LPV}. The model matrices $A,B,C$ are not further optimized after LTI identification, and the remaining parameters to tune are 311 for the deterministic plant model, trained by setting $\rho_\mathrm{w}=2\cdot10^{-3}$. After training the plant model, this is frozen and the remaining 8 parameters of the noise model are trained to get the combined plant+noise model. The results obtained by training the models from 500 initial guesses, with $\rho_\mathrm{w}=2\cdot10^{-4}$,
are summarized in Table~\ref{tab:NL-gyroscope}. We can see that there is a slight improvement on training data with respect to the LPV self-scheduled structure when only estimating a process model, and an improvement on both training and test data in the case a combined model (\combined) is used.

\begin{table}[h!]
\begin{center}
\caption{Identification results for the CMG system (Section \ref{benchmark:cmg})  with NL models.}
\label{tab:NL-gyroscope}
\begin{tabular}{@{}l|c|c|r|r|c|r@{}}
  & $\hat{n}_\mathrm{x}$ & $\hat{n}_\mathrm{z}$ & BFR train. & BFR test & type & time\\
     \hline
  {\tt true} &   5 & -   &  98.16\% &  98.19\% & sim & - \\
\hline
\plantonlyt &
5 & 0 &  96.75\% &  96.12\% & sim &  42.10\sps s\\
\hline
\multirow{2}{*}{\combinedt} & 5 & 2 &  96.66\% &  96.12\% & sim & \multirow{2}{*}{47.78\sps s}\\
& 5 & 2 &  97.82\% &  97.84\% & pred\\
\hline
\end{tabular}
\end{center}
\end{table}

\section{Conclusions} \label{sec:6}
This paper has introduced a novel general identification approach that is based on a theoretically well-founded separate parameterization of the deterministic process and stochastic noise dynamics in the model structure. This not only allows the extension of the estimation concept of the LTI framework to general nonlinear problems, but also to treat the estimation of a whole spectrum of process and noise models from LTI, to LPV and NL dynamics and various noise scenarios encoded by the noise model, by seamlessly tuning the complexity of the model estimates in both of these aspects. 

Thanks to the use of efficient L-BFGS nonlinear optimization and auto-\-dif\-fer\-en\-ti\-ation libraries, the proposed method has been demonstrated to exhibit high computational efficiency in estimating nonlinear dynamics through ANN-based parameterization of the functional relationships, compared to other (deep)-learning-based ANN methods. Besides theoretical guarantees of asymptotic consistency, the introduced approach has shown to provide equivalent and often superior quality of fit than many existing state-of-the art methods for LTI, LPV, and NL identification, with limited needs of hyperparameter tuning. We believe that our proposed work opens the door for a new generation of system identification methods, capable of handling a broad spectrum of model classes with the same simplicity and computational efficiency that users are used to within the LTI framework.      

\section{Acknowledgements}

The authors would like to thank Gerben Beintema and Maarten Schoukens for fruitful discussions on the results. 

\newcommand{\noopsort}[1]{}

\appendix

\section{Proof of Theorem \ref{thm:separation}} \label{proof:separation}
By taking $w_k=[ \ x_k^\top \ \  z_k^\top\ ]^\top$ and $f_\mathrm{w} = [\ f_\mathrm{x}^\top \ \  f_\mathrm{z}^\top\ ]^\top $ with $g_\mathrm{w} = [\ g_\mathrm{x}^\top \ \  g_\mathrm{z}^\top\ ]^\top$, \eqref{eq:system:proc}-\eqref{eq:system:noise} can always be written in the form~\eqref{model:innovation}. For the reverse statement, define a state vector $x_k$ with $n_\mathrm{x}=n_{\mathrm{w}}$ components
and its state-update function $f_\mathrm{x}:\rr^{n_\mathrm{x}}\times\rr^{n_\mathrm{u}}\to\rr^{n_\mathrm{x}}$ as 
\begin{equation}x_{k+1}=f_\mathrm{x}(x_k,u_k)\triangleq f_\mathrm{w}(x_k,u_k,0), \end{equation}
with $x_{0}=w_{0}$. Then, it follows that 
\begin{equation}w_{k+1}-x_{k+1}=f_\mathrm{w}(w_k,u_k,e_k) - f_\mathrm{x}(x_k,u_k) \end{equation}
for all $k\in\mathbb{Z}_0^+$. Now, let us define the state $z_k\triangleq w_{k}-x_{k}$
of the inverse noise model, with $z_0=0$ and
$n_\mathrm{z}=n_{\mathrm{w}}$, and its state-update function $f_\mathrm{z}:\rr^{n_\mathrm{z}}\times\rr^{n_\mathrm{x}}\times\rr^{n_\mathrm{u}}\times\rr^{n_\mathrm{y}}\to\rr^{n_\mathrm{z}}$ as
\begin{align}
z_{k+1}&=f_\mathrm{z} (z_k,x_k,u_k,e_k)\notag \\
&\triangleq f_\mathrm{w}(z_k+x_k,u_k,e_k) - f_\mathrm{w}(x_k,u_k,0).  \label{eq:fz:const}
\end{align}
Similarly, take
\begin{equation}y_{\mathrm{o},k} = g_\mathrm{x}(x_k,u_k) \triangleq g_\mathrm{w}(x_k,u_k) \end{equation}
and define the output function $g_\mathrm{z} (z_k,x_k,u_k) \triangleq g_\mathrm{w}(z_k+x_k,u_k) - g_\mathrm{w}(x_k,u_k)$ such that
\begin{align}
v_k&=y_{k}-y_{\mathrm{o},k} \notag\\
 &=g_\mathrm{w}(z_k+x_k,u_k) +e_k  - g_\mathrm{w}(x_k,u_k) \notag\\
&= g_\mathrm{z} (z_k,x_k,u_k) + e_k.  \label{eq:gz:const}
\end{align}
Note that $f_\mathrm{z}(0,\centerdot,\centerdot,0)\equiv 0$, as by \eqref{eq:fz:const},
\begin{equation}
 f_\mathrm{z} (0,x_k,u_k,0) = f_\mathrm{w}(0+x_k,u_k,0) - f_\mathrm{w}(x_k,u_k,0)=0, 
\end{equation}
and $g_\mathrm{z}(0,\centerdot,\centerdot)\equiv 0$, as by \eqref{eq:gz:const}, $g_\mathrm{w}(0+x_k,u_k)  - g_\mathrm{w}(x_k,u_k)=0$.

This concludes the proof with $n_\mathrm{x}=n_\mathrm{z}=n_\mathrm{w}$. Note that such a state construction can result in a non-minimal state realization for $G_\mathrm{o}$ and $H_\mathrm{o}$. Hence, the corresponding realizations might be further reduced by finding invertible state-transformations with $n_\mathrm{x},n_\mathrm{z}\leq n_\mathrm{w}$ such that the input/output maps corresponding to $G_\mathrm{o}$ and $H_\mathrm{o}$ remain unchanged. 
\hfill $\blacksquare$ 

\end{document}